\documentclass[reqno,12pt]{amsart}

\usepackage{amsmath,amssymb,amsthm,amsfonts}
\usepackage[numbers]{natbib}
\usepackage{libertine}
\usepackage{libertinust1math}
\usepackage[english]{babel}
\usepackage[T1]{fontenc}
\usepackage[utf8]{inputenc} 
\usepackage{enumitem} 
\usepackage[usenames]{xcolor}
\usepackage[babel]{csquotes} 
\usepackage[all]{xypic}
\usepackage{tikz-cd}
\usepackage{ebproof}
\usepackage{accents} 
\usepackage{soul} 
\usepackage[hidelinks]{hyperref}
\usepackage{stmaryrd}
\allowdisplaybreaks

\newtheorem{theorem}{Theorem}[section]

\newtheorem{lemma}[theorem]{Lemma}

\newtheorem{corollary}[theorem]{Corollary}

\theoremstyle{definition}
\newtheorem{example}[theorem]{Example}
\newtheorem{definition}[theorem]{Definition}
\newtheorem{remark}[theorem]{Remark}

\newcommand{\mc}[1]{\mathcal{#1}}

\newcommand{\mb}[1]{\mathbf{#1}}

\newcommand{\mr}[1]{\mathrm{#1}}

\newcommand{\nt}{\Rightarrow}
\newcommand{\Set}{\mb{Set}}

\newcommand{\set}[1]{\left\{#1\right\}}

\newcommand{\cv}{\operatorname{\downarrow}}

\newcommand{\scomp}[2]{ \{\, #1 \mid #2 \,\}}

\newcommand{\op}{^{\operatorname{op}}}

\newcommand{\id}{\operatorname{id}}
\newcommand{\pair}[1]{\langle #1 \rangle}

\newcommand{\eff}{\Leftrightarrow}
\newcommand{\conjt}{\operatorname{\&}}

\newcommand{\gb}{\rightleftarrows}

\newcommand{\ov}[1]{\overline{#1}}

\newcommand{\qsi}[1]{\widetilde{#1}}

\newcommand{\inv}{^{-1}}

\newcommand{\hook}{\hookrightarrow}

\newcommand{\ct}[1]{\underrightarrow{\lim}_{#1}}
\renewcommand{\lim}{\operatorname{Lim}}

\newcommand{\UltL}{\mb{UltL}}
\newcommand{\UltR}{\mb{UltR}}
\newcommand{\Ult}{\mb{Ult}}
\newcommand{\UltSet}{\mb{UltSet}}
\newcommand{\UltPos}{\mb{UltPos}}
\newcommand{\CUltPos}{\mb{CUltPos}}

\newcommand{\UltRPos}{\mb{UltRPos}}
\newcommand{\UltLPos}{\mb{UltLPos}}
\newcommand{\CUltLPos}{\mb{CUltLPos}}

\newcommand{\idl}{\mathrm{Idl}}
\newcommand{\fil}{\mathrm{Fil}}
\newcommand{\DL}{\mathbf{DL}}

\newcommand{\Top}{\mathbf{Top}}

\newcommand{\Mod}{\mathrm{Mod}}

\newcommand{\Pretopos}{\mathbf{Pretopos}}

\newcommand{\Pos}{\mathbf{Pos}}
\newcommand{\CPos}{\mathbf{CPos}}

\newcommand{\Loc}{\mathbf{Loc}}
\newcommand{\CLoc}{\mathbf{CLoc}}
\newcommand{\Bool}{\mathbf{Bool}}

\newcommand{\Comp}{\mb{Comp}}
\newcommand{\Stone}{\mb{Stone}}
\newcommand{\Pries}{\mb{Pries}}

\newcommand{\CTop}{\mb{CTop}}

\newcommand{\clc}[1]{\mathrm{Cl}_{\downarrow}(#1)}
\newcommand{\cld}[1]{\mathrm{Cl}_{\uparrow}(#1)}
\newcommand{\clcd}[1]{\mathrm{Cl}_{\downarrow}^{\uparrow}(#1)}
\newcommand{\cl}[1]{\mathrm{Cl}(#1)}
\newcommand{\dv}{\operatorname{\uparrow}}
\newcommand{\pt}{\mr{pt}}
\newcommand{\elem}[2]{\pair{ #1 , #2 }}

\title{Ultraposet, Distributive Lattice, \\ and Coherent Locale}


\author{Lingyuan Ye}

\date{\today}
  
\begin{document}

\begin{abstract}
  In this paper, we provide an alternative description of the duality result for distributive lattices and coherent locales using \emph{ultraposet}. In particular, we show that there are fully faithful embeddings from the opposite of the category of distributive lattices into the category of ultraposets with ultrafunctors, and from the category of coherent locales into the category of ultraposets with left ultrafunctors. We also define the notion of \emph{zero-dimensional} ultraposets, which characterises the essential image of these embeddings. 
\end{abstract}

\maketitle              

\section{Introduction}\label{sec:intro}

Makkai's 1987 work~\cite{makkai1987stone} on the duality of first-order logic is an essential result for categorical logic. It proves the following socalled \emph{conceptual completeness}:

\begin{theorem}[Conceptual Completeness for First-Order Logic]\label{thm:conceptual}
  The 2-category $\Pretopos$ of (small) pretoposes fully faithfully embeds into the 2-category $\Ult$ of \emph{ultracategories} and \emph{ultrafunctors},
  \[ \Pretopos\op \hook \Ult. \]
\end{theorem}

Every first-order logic have an associated \emph{syntactic category} (cf. e.g.~\cite{caramello2018theories}), and through this construction we may view $\Pretopos$ as the category of first-order theories, up to elimination of imaginaries (c.f.~\cite{harnik2011model}). The embedding from $\Pretopos\op$ to $\Ult$ is then given by taking each theory/pretopos $T$ to its category of \emph{models} $\Mod(T)$.

By name, ultracategories are categories equipped with an ultraproduct structure, and what makes $\Mod(T)$ an ultracategory is exactly the following fundamental theorem in model theory:

\begin{theorem}[\L os Ultraproduct Theorem]\label{thm:los}
  If $\set{M_s}_{s\in S}$ is a family of $T$-models indexed by a set $S$, and $\mu$ is an ultrafilter on $S$, then the ultraproduct of $\set{M_s}_{s\in S}$ w.r.t. $\mu$ is again a $T$-model.
\end{theorem}

From a categorical point of view, the ultraproduct of the family $\set{M_s}_{s\in S}$ w.r.t. $\mu$ is given by the \emph{filtered colimit} $\ct{A\in\mu}\prod_{s\in A}M_s$. The proof of Theorem~\Ref{thm:los} can be found in any standard model theory textbook, e.g.~\cite{hodges1997shorter}. 

Recently, Lurie~\cite{lurie2018ultracategories} has provided an alternative axiomatisation of the notion of ultracategories, \emph{different} from that of Makkai's, and has extended the result to coheret toposes as well,
\[
  \xymatrix{
    \Pretopos\op \ar@{>->}[d] \ar@{^{(}->}[r] & \Ult \ar@{>->}[d] \\
    \mb{CTopos} \ar@{^{(}->}[r] & \UltL
  }
\]
$\mb{CTopos}$ is the 2-category of coherent toposes and geometric morphisms. $\UltL$ has the same objects as $\Ult$, viz. ultracategories, but the morphisms consist of what Lurie calls \emph{left} ultrafunctors. The notion of left ultrafunctors is only present in Lurie's work, and it corresponds certain \emph{lax} version of ultrafunctors, under Lurie's axiomatisation. There is a fully faithful embedding of the 2-category $\mb{CTopos}$ into $\UltL$, by taking each coheret topos $\mc E$ to its category of points $\pt(\mc E)$.

There is a dual \emph{oplax} version of ultrafunctors as well, and they are denoted by Lurie as \emph{right} ultrafunctors. The 2-category $\UltR$ also plays an important role in Lurie's treatment of ultracategories.

The vertical arrows in the above diagramme are faithful, but \emph{not} full. Not every geometric morphism between coherent toposes are induced by a coherent functor between their subcategory of compact objects, and not every ultrafunctor is a left ultrafunctor. Thus, from a topos theoretic point of view, left ultrafunctors are in some sense more fundamental than ultrafunctors, and this is one of the main benefits of Lurie's reaxiomatisation of ultracategories.

In this paper, we perform a \emph{de}categorification of the above results. We will follow Lurie's notion of ultracategories, and (left/right) ultrafunctors, but we specifically focus on \emph{ultraposets}, which are ultracategories with their underlying categories being posets. Our main goal in this paper is to describe a similar duality between distributive lattices, coherent locales, and ultraposets: The opposite of the 2-category $\DL$ of distributive lattices fully faithfully embeds into the 2-category $\UltPos$ of ultraposets and ultrafunctors, and the 2-category $\CLoc$ of coherent locales fully faithfully embeds into the 2-category $\UltLPos$ of ultraposets and \emph{left} ultrafunctors.

The classical theorem of Stone duality also fits into this picture. The category $\Bool$ of Boolean algebras is a full subcategory of $\DL$, and its image under the embedding lies in \emph{ultrasets}, which are \emph{discrete} ultraposets. In~\cite{lurie2018ultracategories}, it is shown that $\UltSet$ is isomorphic to the category $\Comp$ of compact Hausdorff spaces, which contains $\Stone$, the category of Stone spaces, as a full subcategory. Thus, we obtain a diagramme as follows,
\[
  \xymatrix{
    \Bool\op \ar@{^{(}->}[d] \ar@{^{(}->}[r] & \UltSet \ar@{^{(}->}[d] \\
    \DL\op \ar@{>->}[d] \ar@{^{(}->}[r] & \UltPos \ar@{>->}[d] \\
    \CLoc \ar@{^{(}->}[r] & \UltLPos
  }
\]

We will also characterise the essential images of these embeddings by intrinsic \emph{topological} properties of ultraposets. These special ultraposets will be called \emph{zero-dimensional}, and when applied to the discrete case, zero-dimensional ultrasets are exactly the zero-dimensional compact Hausdorff spaces in the topological sense, viz. Stone spaces. This way, we obtain a new duality discritption for distributive lattices using zero-dimensional ultraposets, and also derive the classical theorem of Stone duality as a special case.

Our long term goal is to study the theory of ultraposets, and ultimately ultracategories, in more detail, to apply to the study of model theory. The full axiomatisation of ultracategories in~\cite{lurie2018ultracategories} is quite complex, but it simplifies dramatically in the case of ultraposets. We hope to gain more insight into the full theory of ultracategories by considering this easier case as a first step.

The structure of this paper is as follows. Section~\Ref{sec:basic} will provide the basic definition of ultraposets and (left/right) ultrafunctors, and discuss some typical examples of ultraposets. In Section~\Ref{sec:top}, we study various topological structures naturally induced by ultraposets, which will be fundamental for the duality result later. We construct the embedding functor from $\DL\op$ to $\UltPos$ and $\CLoc$ to $\UltLPos$ in Section~\Ref{sec:uldiscoh}, and show that both distributive lattices and coherent locales can be reconstructed from their correponding ultraposets in Section~\Ref{sec:reconstruct}. Section~\Ref{sec:dual} then proves that these embeddings are fully faithful. We provide an intrinsic characterisation of the essential images of these functors in Section~\Ref{sec:sep}, thus complete the description of the duality theorems. Finally, in Section~\Ref{sec:concfuture} we provide a summary, and discuss various future directions. 

\section{Basic Theory of Ultraposets}\label{sec:basic}

In this section, we define the notion of ultraposets and (left/right) ultrafunctors between them, and see they naturally form three 2-categories $\UltPos$, $\UltLPos$, and $\UltRPos$. As we've mentioned, our treatment of ultraposets will essentially follow~\cite{lurie2018ultracategories}. However, we adopt slightly different notation choices for the ultrastructure, for better readability.

We begin by recalling the basics of ultrafilters. An ultrafilter $\mu$ on a set $S$ is an upward closed subset of $\wp(S)$, which contains $S$, is closed under finite intersection, and for any two complemented subset of $S$, exactly one of them is in $\mu$.

Recall that we have a functor $\beta : \Set \to \Set$, sending each $S$ to its set of ultrafilters $\beta S$ on $S$. For any function $i : S \to T$, we will denote the map from $\beta S$ to $\beta T$ as
\[ i_{*} : \beta S \to \beta T, \]
where for any $\mu\in\beta S$, we have
\[ i_{*}\mu = \scomp{V\subseteq T}{i\inv(V) \in \mu}. \]
We will denote $i_{*}\mu$ as the \emph{pushforward} of $\mu$ along $i$. Notice that if $i : S \hook T$ is injective, then $i_{*} : \beta S \hook \beta T$ is also injective, and we can identify $\beta S$ as the subspace of $\beta T$ consisting of those $\nu \in \beta T$ that $i(S) \in \nu$. 

There is actually a unique \emph{monad} structure on $\beta$, with $\delta : \id \nt \beta$ and $\gamma : \beta\circ\beta \nt \beta$.\footnote{Uniqueness of the monad structure follows from the fact that $\beta$ is \emph{terminal} among functors from $\Set$ to $\Set$ that preserves finite coproducts. This is first proven in~\cite{borger1987coproducts}.} Concretly, given any set $S$,
\begin{itemize}
\item For any $s\in S$, $\delta_s = \scomp{A\subseteq S}{s\in A}$.
\item For any ultrafilter $\theta$ on $\beta S$, $\gamma(\theta) = \scomp{A\subseteq S}{\scomp{\mu\in\beta S}{A \in \mu}\in\theta}$. 
\end{itemize}
The image of $\delta$ are denoted as \emph{principal} ultrafilters.

Since each $\beta S$ is a free $\beta$-algebra, given any $\mu : T \to \beta S$, there is an induced morphism from $\beta T \to \beta S$, and we write this map as follows,
\[ \elem\mu{-} : \beta T \to \beta S. \]
Concretely, for any $\nu\in\beta T$, the ultrafilter $\pair{\mu,\nu}$ on $\beta S$ is given by
\[ A \in \pair{\mu,\nu} \eff \mu\inv(\beta A) = \scomp{t}{A\in\mu_t} \in \nu. \]
We slightly extends it to a map of the following type,
\[ \elem\mu{-} : \beta T^R \to \beta S^R. \]
where for any function $\nu : R \to \beta T$, we have
\[ \elem\mu\nu = r \mapsto \elem\mu{\nu_r}. \]
The pairing $\elem{-}{-}$ is essentially the Kleisli structure of the monad. Under this notation, the Kleisli identities imply that for any $\mu : T \to \beta S$,
\[ \elem\mu\delta = \mu, \]
and furthermore for any $\mu : T \to \beta S$, $\nu : R \to \beta T$, and $\lambda : W \to \beta R$,
\[ \elem\mu{\elem\nu\lambda} = \elem{\elem\mu\nu}\lambda. \]
Under this notation, the action of the functor $\beta$ on morphism can also be expressed as follows,
\[ i_{*}\nu = \elem{\delta i}{\nu}. \]

\begin{remark}\label{rem:converge}
  From a topological perspective, each $\beta S$ is a compact Hausdorff space, and it is in fact the Stone-\u{C}ech compatification of the discrete space $S$. It is well-known that for any compact Hausdorff space $X$, any ultrafilter $\mu$ on $X$ has a \emph{unique} convergence point in $X$, viz. a point $x\in X$ where $\tau_x \subseteq \mu$, with $\tau_x$ being the set of all open neighbourhoods of $x$. This map $\beta X \to X$ actually makes $X$ an \emph{algebra} for the monad $\beta$, and in fact, it is shown already in the work of Ernest Manes~\cite{Manes1976} that there is an isomorphism of categories,
  \[ \Set^{\beta} \cong \Comp. \]
  We refer the readers to~\cite[Ch. 3]{goldbring2022ultrafilters} for a more detailed account of the connection between ultrafilter and topology. The pairing $\pair{\mu,\nu}$ henceforth can be seen as the unique convergence point in $\beta S$ of the ultrafilter $\mu_{*}\nu$. Suggestively,~\cite{lurie2018ultracategories} uses the notation $\int\mu d\nu$ to denote this convergence point. For better readability, we have adopted the pairing notation instead of integration in this paper.
\end{remark}

We now tend to the definition of ultraposets. An ultraposet $P$ will be a poset $P$, equipped with similar operations that are available on $\beta S$ as described above, and all $\beta S$, viewed as a discrete poset, will indeed be instances of ultraposets:

\begin{definition}[Ultraposet]\label{def:ultraposet}
  An \emph{ultraposet} is a poset $P$ equipped with the following data: For any set $S$ and $T$, there is a \emph{monotone} pairing function providing the ultraproduct structure on $P$, 
  \[ \elem-- : P^S \times \beta S^T \to P^T, \]
  where $P^S$ and $P^T$ has the point-wise order, and $\beta S^T$ is discrete. We require this pairing to be point-wise on the second entry, in the sense that for any $f : S \to P$, $\mu : T \to \beta S$ and any $i : W \to T$,
  \[ \pair{f,\mu} \circ i = \pair{f,\mu\circ i}. \]
  The pairing function should furthermore satisfy the following properties: For any $f : S \to P$ and any $\mu : T \to \beta S$, $\nu : R \to \beta T$,
  \begin{enumerate}
  \item \label{un} \emph{Unity}: The principal ultrafilers are units,
    \[ \elem f\delta = f. \]
  \item \label{lt} \emph{Lax Associativity}: The pairing is associative upto an inequality,
    \[ \pair{f,\pair{\mu,\nu}} \le \pair{\pair{f,\mu},\nu}. \]
  \item \label{lc} \emph{Locality}: If $i : W \hook S$ is injective, then
    \[ \pair{f,i_{*}\nu} = \pair{f,\pair{\delta i,\nu}} = \pair{\pair{f,\delta i},\nu} = \pair{fi,\nu}. \]
  \end{enumerate}
\end{definition}

We provide some intuitive explanation on the definition of an ultraposet. Firstly, the pairing being point-wise on the second entry just reflects our \emph{notation} choice, which has nothing essential to do with the structure of ultraposets. Alternatively, we can also define the pairing function to be a map of type
\[ \pair{-,-} : P^S \times \beta S \to P, \]
and this is indeed the choice in~\cite{lurie2018ultracategories}. Any such pairing function can be extended point-wise to one with type $P^S \times \beta S^T \to P^T$, by setting 
\[ \pair{f,\mu} = t \mapsto \pair{f,\mu_t}, \]
for any $f : S \to P$ and any $\mu : T \to \beta S$. We find our notation more succinct and clear, especially when expressing and using lax associativity.

The remaining three properties are at the heart of ultraposets:
\begin{itemize}
\item The unity condition is best understood from the perspective of convergence (cf. Remark~\Ref{rem:converge}). Let $\delta_s$ be the principal ultrafilter for $s\in S$. The family $\set{f_s}_{s\in S}$ intuitively converge to the point $f_s$ under the ultrafilter $\delta_s$, because it contains the singleton subset $\set{s}$.
\item The lax associativity condition would become more sensible once we have seen more examples of ultracategories. For now, notice that if the poset is \emph{discrete}, then lax associativity is actually equivalent to strict associativity, and in fact the pairing for $\beta S$, or more generally for any compact Hausdorff space, is strictly associative.
\item The locality condition is also crucial. It asserts that the pairing operation $\pair{f,\mu}$ is \emph{local}. We have mentioned that, when $i : W \hook S$ is injective, we can view $\beta W$ as a subspace of $\beta S$, under the inclusion $i_{*} : \beta W \hook \beta S$. Given any $\mu\in\beta S$, $\mu$ lies in the image of $i_{*}$, viz. $\mu = i_{*}\nu$ for some $\nu$, iff $W \in \mu$. When this is the case, the convergence point $\pair{f,\mu}$ should be equal to the convergence point $\pair{fi,\nu}$ when we restrict the family on $S$ to the subfamily on $W$, and restrict the ultrafilter $\mu$ on $\nu$.
\end{itemize}

Notice that for any ultraposet $P$, given $f : S \to P$ and an arbitrary function $g : W \to S$, if $\lambda\in\beta W$, then by lax associativity,
\[ \pair{f,g_{*}\mu} = \pair{f,\pair{\delta g,\mu}} \le \pair{\pair{f,\delta g},\mu} = \pair{fg,\mu}. \]
It follows that the operators $- \circ g$ of pre-composition with $g$, and $g_{*}(-)$ of pushforwarding an ultrafilter, act like a pair of adjoint operators for the pairing function $\pair{-,-}$. However, for arbitrary $g$, this is not quite the case, because we have an inequality; only when $g$ is injective it becomes a true equality.

For ultraposets, there are naturally three different notions of morphisms:

\begin{definition}[Ultrafunctors]\label{def:ultrafunc}
  Given two ultraposet $P,Q$, a \emph{left ultrafunctor} from $P$ to $Q$ is a monotone function $\varphi : P \to Q$, such that for any set $S$, any family $f : S \to P$ and any $\mu : T \to \beta S$,
  \[ \varphi\pair{f,\mu} \le \pair{\varphi f,\mu}. \]
  Here the order is point-wise on $P^T$. Dually, a \emph{right ultrafunctor} is a functor $\varphi : P \to Q$, satisfying
  \[ \pair{\varphi f,\mu} \le \varphi\pair{f,\mu}, \]
  $\varphi$ is an \emph{ultrafunctor} iff it is both a left and right ultrafunctor, i.e. iff the above inequalities are in fact equalities.
\end{definition}

We will use $\UltLPos$ (resp. $\UltRPos$) to denote the category of ultraposets and left (resp. right) ultrafunctors. They have a common wide subcategory $\UltPos$ of ultraposets and ultrafunctors. For any two posets, the set of monotone maps between them has a natural point-wise order, thus both $\UltLPos$, $\UltRPos$ and $\UltPos$ will be $\Pos$-enriched, where $\Pos$ is the category of posets. Hence, they can be naturally considered as 2-categories.

\begin{remark}\label{rem:compareuc}
  Our definitions of ultraposets and (left/right) ultrafunctors are indeed special cases of the more general notions of ultracategories and (left/right) ultrafunctors defined in~\cite{lurie2018ultracategories}, restricted to posets. However, Definition~\Ref{def:ultraposet} and~\Ref{def:ultrafunc} are dramatically simpler than the more general case. This is mainly because, when the category is not a poset, there will be further non-trivial \emph{coherence} conditions on the structure of ultraproducts in an ultracategory. However, since a poset is \emph{skeletal} and \emph{thin}, i.e. isomorphic objects are actually equal and between two objects there are at most one morphism, all the coherence conditions will be trivial.
\end{remark}

We will now discuss some typical examples of ultraposets. As mentioned before, discrete ultraposets are exactly compact Hausdorff spaces:

\begin{example}[Ultrasets]
  If $P$ is a \emph{discrete} ultraposet, i.e. the order on $P$ is the equality relation, then we say $P$ is an \emph{ultraset}. For an ultraset, every inequality by definition coincide with strict equality. In particular, the notion of left/right ultrafunctors coincide with the notion of ultrafunctors in this case, and we denote the category of ultrasets as $\UltSet$. In fact, Lurie has shown that there is an isomorphisms of categories,
  \[ \UltSet \cong \Comp. \]
  The construction making a compact Hausdorff space $X$ into an ultraset is essentially provided by Remark~\Ref{rem:converge}. On the other hand, for any ultraposet, and ultraset in particular, there is a notion of \emph{closed} subsets, as we will see in Section~\Ref{sec:top}. The closed subsets of $X$ in terms of its ultrastructure are exactly the closed subsets of $X$ from its topology. For more details, see~\cite{lurie2018ultracategories}.
\end{example}

In particular, each $\beta S$ will be an ultraset. Recall the lax associativity condition in Definition~\Ref{def:ultraposet} of ultraposet. For any $f : S \to P$, the pairing on $P$ induces a function
\[ \varphi_f := \pair{f,-} : \beta S \to P, \]
and lax associativity exactly asserts that this function $\varphi_f$ is a \emph{left ultrafunctor} for any $f$: For any $\mu : T \to \beta S$ and any $\nu\in\beta T$, $\varphi_f$ being a left ultrafunctor exactly means
\[ \pair{f,\pair{\mu,\nu}} = \varphi_f\pair{\mu,\nu} \le \pair{\varphi_f\mu,\nu} = \pair{\pair{f,\mu},\nu}. \]
This is a first sign that the notion of left ultrafunctors might be more fundamental than the notion of ultrafunctors.

More generally, we can equip any compact Hausdorff space with a \emph{closed} partial order to make it into an ultraposet:

\begin{lemma}\label{lem:compord}
  Let $X$ be a compact Hausdorff space viewed as an ultraset. If $\le$ is a partial order on $X$ which is \emph{closed} in the space $X \times X$, $(X,\le)$ will be an ultraposet.
\end{lemma}
\begin{proof}
  The ultrastructure on $X$ evidently satisfies Conditions~(\Ref{un}),~(\Ref{lt}) and~(\Ref{lc}) in Definition~\Ref{def:ultraposet} as well for the poset $(X,\le)$, basically because $\le$ is reflexive. Hence, we only need to verify the pairing function is also monotone. Let $f,g : S \to X$ be two functions that $f \le g$ in the point-wise order on $X^S$ when $X$ is equipped with $\le$. For any $\mu\in\beta S$, let $\pair{f,\mu}$ be $x$, and $\pair{g,\mu}$ be $y$. By Remark~\Ref{rem:converge}, $x$ is the unique point satisfying
  \[ \forall U\in\tau_x, f\inv(U) \in \mu, \]
  and similarly for $g$. Since $\le$ is closed, if $x \not\le y$, there will be open neighbourhoods $U_x\in\tau_x$ and $U_y\in\tau_y$ that
  \[ (U_x \times U_y) \cap \operatorname{\le} = \emptyset. \]
  This way, it follows that $f\inv(U_x) \in \mu$ and $g\inv(U_y) \in \mu$, thus
  \[ f\inv(U_x) \cap g\inv(U_y) \in \mu. \]
  In particular, $f\inv(U_x) \cap g\inv(U_y) \neq \emptyset$, thus we can find $s \in S$ in this intersection. Then by our construction,
  \[ f(s) \in U_x \conjt g(s) \in U_y \nt f(s) \neq g(s), \]
  contradicting our assumption that $f \le g$. Thus, we must have
  \[ \pair{f,\mu} \le \pair{g,\mu}, \]
  and this completes the proof that $(X,\le)$ equipped with this ultrastructure is an ultraposet.
\end{proof}

The notion of a compact Hausdorff space equipped with a closed partial order is already studied in Nachbin's work~\cite{nachbin1965topology}. In modern terminology it is usually called \emph{compact ordered space}, or \emph{compact pospace}. There is a natural (2-)category $\mb{CompOrd}$ of compact ordered spaces, with morphisms being continuous and monotone maps, with point-wise order on morphisms. From Lemma~\ref{lem:compord} and the isomorphism $\Comp \cong \UltSet$, it follows that there is a fully faithful embedding of (2-)categories as follows,
\[ \mb{CompOrd} \hook \UltPos. \]
We will come back to these examples in Section~\Ref{sec:sep}.

Another family of examples comes from complete lattices. Recall from Section~\Ref{sec:intro} that the ultraproduct of a family of models $\set{M_s}_{s\in S}$ w.r.t. $\mu\in\beta S$ can be viewed as the filtered colimit $\ct{A\in\mu}\prod_{s\in A}M_s$. We can mimic this definition, and equip any complete lattice with an ultrastructure:

\begin{example}[Canonical Ultrastructure]
  Every complete lattice $P$ can be made into an ultraposet. Given any $f : S \to P$ and any $A \subseteq S$, we use $f^A$ to denote the meet $\bigwedge_{s\in A}fs$ in $P$. For any $\mu\in\beta S$, we define
  \[ \pair{f,\mu} := \bigvee_{A\in\mu}f^A. \]
  Notice that the colimit is actually indexed by the direct poset $\mu\op$, if we view the order in $\mu$ as given by inclusion of subsets. This is because if $A \subseteq B$, then by definition $f^B = \bigwedge_{s\in B}fs \le \bigwedge_{s\in A}fs = f^A$. From the above definition, it is evident that the pairing $\pair{-,\mu}$ will be monotone. We verify in detail that it also satisfies the remaining three conditions:
  \begin{itemize}
  \item Unity is obviously satisfied, because for any $s\in S$, $\set{s}$ will be the bottom element in $\delta_s$, thus a top element in $\delta_s\op$. This implies the join is simply equal to $f^{\set s} = fs$.
  \item For lax associativity, given $\mu : T \to \beta S$ and $\nu\in\beta T$, for any $A\subseteq S$, by definition we have
    \[ A \in \pair{\mu,\nu} \eff \mu\inv(\beta A) = \scomp{t}{A\in\mu_{t}} \in \nu. \]
    This way, for any $A \in \pair{\mu,\nu}$ and any $t\in\mu\inv(\beta A)$, $A\in\mu_t$, and thus by definition of $\pair{f,\mu_t}$ we have
    \[ f^A \le \pair{f,\mu_t}. \]
    It follows that
    \[ f^A \le \bigwedge_{t\in\mu\inv(\beta A)}\pair{f,\mu_t} = \pair{f,\mu}^{\mu\inv(\beta A)}. \]
    This means for any $A\in\pair{\mu,\nu}$, we can find some subset $\mu\inv(\beta A) \in \nu$, that $f^A$ is smaller than $\pair{f,\mu}^{\mu\inv(\beta A)}$. Hence, by definition we obtain
    \[ \pair{f,\pair{\mu,\nu}} \le \pair{\pair{f,\mu},\nu}. \]
  \item Locality is also easy, because if $i : T \hook S$ is an inclusion, then $\nu\op$ will be \emph{final} in the $i_{*}\nu\op$, hence they have the same join.
  \end{itemize}
\end{example}

The above will be denoted as the \emph{canonical}, or \emph{categorical}, ultrastructure on a complete lattice. Notice that in our proof of lax associativity, \emph{a priori} there is no reason that the inequality $\pair{f,\pair{\mu,\nu}} \le \pair{\pair{f,\mu},\nu}$ is actually an equality, thus we do need this generality in the definition of ultraposets if we want to include these examples.

One most important instance of an ultraposet for this paper is the ultraposet $\mb 2 = \set{0<1}$. On one hand, it is a complete lattice, thus can be equipped with the canonical ultrastructure. For any family $f : S \to \mb 2$ and any $\mu\in\beta S$, we observe that
\[ \pair{f,\mu} = \bigvee_{A\in\mu}f^A = 1 \eff \exists A\in\mu,\,A \subseteq f\inv(1) \eff f\inv(1) \in \mu. \]
Since $f\inv(1)$ and $f\inv(0)$ are disjoint and they cover $S$, it follows that exactly one of them will be in $\mu$, and the above characterisation suggests that $\pair{f,\mu} = i$ iff $f\inv(i) \in \mu$, for $i \in \set{0,1}$.

What we may also observe is that, $\mb 2$ as a set is also a compact Hausdorff space with the discrete topology, and its canonical ultrastructure actually \emph{coincide} with the ultrastructure provided by this topology. This implies that $\mb 2$ is also an instance from $\mb{CompOrd}$, and indeed $\le$ is closed in $\mb 2$, because it is discrete. The later part of this paper will depend extensively on the ultrastructure on $\mb 2$.

The canonical ultrastructure also extends to a fully faithful embedding of a suitable (2-)category into the (2-)category $\UltL$:

\begin{lemma}\label{lem:leftscott}
  Given a monotone map $\varphi : P \to Q$ between two complete lattices, view $P,Q$ as ultraposets equipped with the canonical ultrastructure, then $\varphi$ is a left ultrafunctor iff it preserves arbitrary directed joins.
\end{lemma}
\begin{proof}
  Let $f : S \to P$ be any function, and let $\mu$ be an element in $\beta S$. If $\varphi$ preserves directed joins,
  \[ \varphi\pair{f,\mu} = \varphi\left(\bigvee_{A\in\mu}f^A\right) = \bigvee_{A\in\mu}\varphi(f^A) \le \bigvee_{A\in\mu}(\varphi f)^A = \pair{\varphi f,\mu}. \]
  It follows that $\varphi$ is a left ultrafunctor. On the other hand, consider a family of objects $\set{p_i}_{i\in I}$ in $P$, indexed by some directed poset $I$. Let $\dv i$ denote the following subset of $I$,
  \[ \dv i := \scomp{j\in I}{i \le j}. \]
  Observe that the family of subsets $\set{\dv i}_{i\in I}$ has the finite intersection property because $I$ is directed. Thus, we can find some ultrafilter $\mu$ on $I$ that extends this family.\footnote{This is a well-known property of ultrafilters. For a proof, see Appendix~\Ref{app}, Lemma~\Ref{lem:maximalideal}.} Then by definition,
\[ \pair{p,\mu} = \bigvee_{A\in\mu}p^A. \]
Notice that since $\dv i\in\mu$ for each $i\in I$, evidently we have
\[ p_i = p^{\dv i} \le \bigvee_{A\in\mu}p^A = \pair{p,\mu}, \]
which implies $\bigvee_{i\in I}p_i\le \pair{p,\mu}$. On the other hand, for any $A\in\mu$,
\[ p^A = \bigwedge_{i\in A}p_i \le \bigvee_{i\in I}p_i. \]
Thus, it follows that
\[ \pair{p,\mu} = \bigvee_{A\in\mu}p^A \le \bigvee_{i\in I}p_i. \]
This proves that $\pair{p,\mu}$ and $\bigvee_{i\in I}p_i$ coincide. Similarly, we also have
\[ \pair{\varphi p,\mu} = \bigvee_{i\in I}\varphi p_I, \]
for any monotone function $\varphi$. Now if $\varphi$ is a left ultrafunctor, we must have
\[ \varphi\left(\bigvee_{i\in I}p_i\right) = \varphi\pair{p,\mu} \le \pair{\varphi p,\mu} = \bigvee_{i\in I}\varphi p_i, \]
which actually implies $\varphi$ preserves this directed join $\bigvee_{i\in I}p_i$. Hence, left ultrafunctors between complete lattices with canonical ultrastructure coincide with Scott-continuous maps between them.
\end{proof}


Let $\CPos$ be the (2-)category of complete lattices with Scott-continuous functions, viz. maps preserving directed joins. Lemma~\ref{lem:leftscott} then implies that there is a \emph{fully faithful} embedding
\[ \CPos \hook \UltLPos. \]
$\CPos$ is a full subcategory of DCPOs, the category of directly complete partial orders with Scott-continuous functions between them. In particular, $\CPos$ contains as a full subcategory the category of \emph{continuous lattices}, which is Cartesian closed, and used by Dana Scott~\cite{scott1976data} to build models of untyped $\lambda$-calculus. See~\cite{abramsky1995domain} for more application of DCPOs in domain theory. There are potentially much more connection between the theory of ultraposets and domain theory, and we will come back to this in Section~\Ref{sec:concfuture}. 

Given posets $P,Q$, let $[P,Q]$ be the poset of all monotone maps from $P$ to $Q$, with point-wise order. Notice that if $P,Q$ are in $\CPos$, $\CPos(P,Q)$ is actually a complete lattice, with joins computed point-wise in $[P,Q]$. This is because the point-wise join of a family of monotone maps that preserves direct joins again preserves directly joins. Then Lemma~\Ref{lem:leftscott} suggests that if we view $P,Q$ as ultraposets, the poset $\UltLPos(P,Q)$ will again be complete, and joins of left ultrafunctors are also computed point-wise there. This is generally true in a more general context:

\begin{lemma}\label{lem:joinpoint}
  For any ultraposets $P,Q$ with $Q$ a complete lattice\footnote{Here we do \emph{not} require $Q$ to be equipped with the canonical ultrastructure.}, the poset $\UltLPos(P,Q)$ (resp. $\UltRPos(P,Q)$) is complete, and its joins (resp. meets) are taken in $[P,Q]$.
\end{lemma}
\begin{proof}
  Suppose we are given a family $\set{\varphi}_{i\in I}$ of left ultrafunctors from $P$ to $Q$, and consider $\varphi$ as their point-wise join. Given any family $f : S \to P$ and any $\mu\in\beta S$, we have
  \[ \varphi\pair{f,\mu} = \bigvee_{i\in I}\varphi_i\pair{f,\mu} \le \bigvee_{i\in I}\pair{\varphi_i f,\mu} \le \pair{\bigvee_{i\in I}\varphi_i f,\mu} = \pair{\varphi f,\mu}. \]
  The first inequality uses the fact that each $\varphi_i$ is a left ultrafunctor, and the second inequality is due to the monotonicity of $\pair{-,\mu}$. Thus, $\varphi$ is a left ultrafunctor, and it is the join of $\set{\varphi_i}_{i\in I}$ in $\UltLPos(P,Q)$. The proof for right ultraposets is completely symmetric.
\end{proof}

Finally, we show how to construct new ultraposets from olds ones:

\begin{example}\label{exm:catconsultrapos}
  We describe products, finite coproducts, and exponentials of ultraposets:
  \begin{itemize}
  \item For any family $\set{P_i}_{i\in I}$, the \emph{product ultraposet} $\prod_{i\in I}P_i$ is simply their product in $\Pos$, equipped with point-wise ultraproduct structure: For any $(f_i)_{i\in I} : S \to \prod_{i\in I}P_i$ and $\mu\in\beta S$, we have
    \[ \pair{(f_i)_{i\in I},\mu} := (\pair{f_i,\mu})_{i\in I}. \]
    It is easy to check that $\prod_{i\in I}P_i$ becomes an ultraposet with the above ultrastructure. In particular, the \emph{terminal} ultraposet is the singleton set $\mb 1$ with trivial order and trivial ultrastructure.

    Furthermore, the point-wise ultrastructure implies that the projection maps are all \emph{ultrafunctors},
    \[ \pi_i : \prod_{i\in I}P_i \to P_i. \]
    It is then straight forward to check that $\prod_{i\in I}P_i$ is indeed the categorical product of $\set{P_i}_{i\in I}$, in both $\UltLPos$, $\UltRPos$, and $\UltPos$.\vspace{1ex}
  \item The \emph{initial} ultraposet is the empty poset $\emptyset$. Since there are no maps from a non-empty set into $\emptyset$, and there are no ultrafilters on $\emptyset$ itself, there is a unique ultrastructure on $\emptyset$.

    Given two ultraposets $P,Q$, their \emph{coproduct} $P\sqcup Q$ is simply the disjoint union of $P$ and $Q$ in $\Pos$. For any map $f : S \to P \sqcup Q$ and any $\mu\in\beta S$, exactly one of the subsets $f\inv(P)$ and $f\inv(Q)$ will be a member of $\mu$, say $f\inv(P) \in \mu$, then we simply assign the ultraproduct $\pair{f,\mu}$ as the ultraproduct $\pair{f|_{f\inv(P)},\mu|_{f\inv(P)}}$ in $P$. In this case, the two inclusions $P,Q \hook P \sqcup Q$ are evidently ultrafunctors as well, and it is also easy to verify that $P \sqcup Q$ is indeed the coproduct of $P,Q$ in $\UltLPos$, $\UltRPos$, and $\UltPos$. \vspace{1ex}
  \item Given any ultraposet $P$ and any poset $Q$, we may equip $[Q,P]$ with a point-wise ultrastructure: For any $f : S \to [Q,P]$ and any $\mu\in\beta S$, we define
    \[ \pair{f,\mu} := q \mapsto \pair{f(q),\mu}. \]
    Verifying $[Q,P]$ is an ultraposet is again straight forward.
  \end{itemize}
\end{example}

In this section, we have defined the 2-categories $\UltLPos$, $\UltRPos$, and $\UltPos$, and discussed some typical examples of ultraposets. The next section discusses the close relationship between ultraposets and certain topological structures, which will be quite useful latter on when we study the duality for distributive lattices and coherent locales.

\section{Topological Structures in Ultraposets}\label{sec:top}

From the isomorphism $\UltSet \cong \Comp$, we should already anticipate a close relationship between ultrastructures and topological structures. In fact, for any ultraposet $P$, the ultraproduct structure on $P$ will induce a compact topology on $P$, by defining the closed sets as follows:

\begin{definition}[Closed Sets of Ultraposets]\label{def:closeset}
  Given an ultraposet $P$, a subset $K \subseteq P$ is \emph{closed} in $P$, iff for any $f : S \to P$, and any $\mu\in\beta S$,
  \[ f\inv(K) \in \mu \nt \pair{f,\mu} \in K. \]
\end{definition}

According to locality of ultraproduct, if $f\inv(K) \in \mu$, then we may simply restrict the ultrafilter $\mu$ to some ultrafilter $\nu$ on $f\inv(K)$, and we must have
\[ \pair{f,\mu} = \pair{f|_{f\inv(K)},\nu}. \]
It follows that $K$ is closed iff for any $f : S \to K \subseteq P$ and any $\mu\in\beta S$, $\pair{f,\mu}$ in $P$ actually lies in $K$.

Evidently, if $K$ is closed in $P$, then viewed as a subposet of $P$, $K$ inherits an ultraposet structure. This applies to ultrafunctors as well:

\begin{lemma}\label{lem:clores}
  Let $\varphi : P \to Q$ be a (left/right) ultrafunctor. Given closed subposet $K \subseteq P$ and $L \subseteq Q$, if $\varphi$ restricts to a morphism from $K$ to $L$, then it will also be a (left/right) ultrafunctor from $K$ to $L$.
\end{lemma}
\begin{proof}
  This is straight forward, because the ultrastructures on $K,L$ are restricted from $P,Q$, respectively.
\end{proof}

The family of closed subsets in any ultraposet actually forms a topolgy:

\begin{lemma}\label{lem:closedtop}
  Closed subsets of any ultraposet $P$ are closed under arbitrary intersection and finite union.
\end{lemma}
\begin{proof}
  Closed under intersection is evident. Given any family $\set{K_i}_{i\in I}$ of closed subsets of $P$, if $f : S \to \bigcap_{i\in I}K_i \subseteq P$, then by each $K_i$ being closed, for any $\mu\in\beta S$ we would have
  \[ \pair{f,\mu} \in K_i, \ \forall i\in I. \]
  It follows that $\pair{f,\mu} \in \bigcap_{i\in I}K_i$, thus the intersection is also closed.

  By definition $P$ itself is closed. Now suppose both $K$ and $L$ are closed subsets in $P$. Consider any family $f : S \to K \cup L \subseteq P$ and any ultrafilter $\mu\in\beta S$. The decomposition $K\cup L$ induces a decomposition $f\inv(K) \cup f\inv(L)$ of $S$, and at least one of them is in $\mu$, say $f\inv(K) \in \mu$. Then by locality of ultraproduct, we must have
  \[ \pair{f,\mu} = \pair{f|_{f\inv(K)},\mu|_{f\inv(K)}}. \]
  Thus by the fact that $K$ is closd, we have $\pair{f,\mu} \in K \subseteq K \cup L$. This implies that $K\cup L$ is closed.
\end{proof}

Given an ultraposet, we will use $\cl P$ to denote the lattice of closed sets of $P$ under inclusion. We furthermore show that the topology $\cl P$ is \emph{compact} for any ultraposet $P$. In the language of closed sets, a topology is compact iff for any family of closed subsets, if each of its finite subfamily has non-empty intersection, then the whole family also has non-empty intersection. Or equivalently, if a family has empty intersection, then it has a finite subfamily that has empty intersection.

\begin{lemma}\label{lem:compact}
  For any ultraposet $P$, $\cl P$ is \emph{compact} considered as closed sets for a topology on $P$.
\end{lemma}
\begin{proof}
  Suppose we have some family $\set{K_i}$ of closed subsets of an ultraposet $P$. Suppose any finite subfamily of $\set{K_i}$ has non-empty intersection. It is well-known that any such family can be extended to an ultrafilter $\mu$ on $P$ that contains $\set{K_i}$.\footnote{For a proof, see Appendix~\Ref{app}, Lemma~\Ref{lem:maximalideal}.} Consider the ultraproduct $\pair{\id,\mu}$. For any $K_i$, since $K_i\in\mu$ and $K_i$ is closed in $P$, by definition we know $\pair{\id,\mu}\in K_i$. This suggests that $\pair{\id,\mu}\in\bigcap_iK_i$, thus the whole family also has non-empty intersection.
\end{proof}

\begin{remark}
  For an ultraset $X$, its closed set $\cl X$ will indeed becomes a compact Hausdorff topology on $X$, and the ultrastructure induced by this topology will coincide with the original ultrastructure on $X$. Again, see~\cite{lurie2018ultracategories} for more details.
\end{remark}

The next natural question to ask is whether the various morphisms between ultraposets are continuous for this topology or not. As we will see, in general neither left nor right ultrafunctors will be continuous for this topology, but ultrafunctors \emph{are}, and they also preserve closed sets:

\begin{lemma}\label{lem:ultracont}
  Any ultrafunctor $\varphi : P \to Q$ is a closed continuous map, when $P,Q$ are equipped with the topology of closed sets $\cl P,\cl Q$, respectively.
\end{lemma}
\begin{proof}
  For any closed set $L$ in $Q$, we show $\varphi\inv(L)$ is also closed. Suppose we are given $f : S \to \varphi\inv(L)\subseteq P$ and $\mu\in\beta S$, by definition $\varphi f : S \to L \subseteq Q$, and thus by $L$ being closed and $\varphi$ being an ultrafunctor,
  \[ \varphi\pair{f,\mu} = \pair{\varphi f,\mu} \in L, \]
  It follows that $\pair{f,\mu}\in\varphi\inv(L)$, hence $\varphi\inv(L)$ is closed.

  We also show $\varphi$ is closed. Let $K \subseteq P$ be a closed subspace. For any map $f : S \to \varphi(K) \subseteq Q$ and any $\mu\in\beta S$, by Choice we can lift $f$ as follows,
  \[
    \xymatrix{
      & S \ar@{-->}[dl]_{\qsi f} \ar[d]^f \\
      P \ar[r]_{\varphi} & Q
    }
  \]
  such that $\qsi f : S \to K \subseteq P$. Then since $K$ is closed and $\varphi$ is an ultrafunctor, it follows that
  \[ \pair{f,\mu} = \pair{\varphi\qsi f,\mu} = \varphi\pair{\qsi f,\mu} \in \varphi(K). \]
  Hence, $\varphi(K)$ is also closed.
\end{proof}

Let $\CTop$ be the category of compact Topological spaces, and let $\CTop_c$ be its wide subcategory where morphisms are restricted to closed continuous maps. Lemma~\Ref{lem:ultracont} then suggests that the construction $\cl-$ extends to a functor as follows,
\[ \cl- : \UltPos \to \CTop_c. \]

From the proof of Lemma~\Ref{lem:ultracont}, we also realise that if $\varphi$ is merely a left or right ultrafunctor, then it may not be continuous for the topology $\cl P$ and $\cl Q$, because the equality bewteen $\varphi\pair{f,\mu}$ and $\pair{\varphi f,\mu}$ may no longer hold, and thus $\varphi\inv (L)$ may not be closed. This is not optimal, because according to our philosophy, left ultrafunctors are more important than ultrafunctors. 

However, the remedy is easy to find. If $\varphi$ is merely a left or right ultrafunctor, we may restrict to those closed subsets that are in addition \emph{downward} or \emph{upward closed}. For any ultraposet $P$, we write $\clc P$ (resp. $\cld P$) for the set of closed and downward (resp. upward) closed subsets of $P$.

What are some basic examples of elements in $\clc P$ and $\cld P$? One might expect that $\cv p \in \clc P$ for any $p\in P$, but this is not the case in general. Consider the constant function $p : S \to P$ at some $p\in P$. For any $\mu\in\beta S$, if we use $p^{\mu}$ to denote $\pair{p,\mu}$, then we have
\[ p = \pair{\id,\delta_p} = \pair{\id,p_{*}\mu} \le p^{\mu}. \]
If the inequality is strict, we would have $p^{\mu} \not\in \cv p$, which implies $\cv p$ might not be closed.

However, for the examples we have considered so far, including ultrasets, compact ordered spaces, and complete lattices with canonical ultrastructures, subsets of the form $\cv p$ will in fact be closed. For an ultraset $X$, for any $p\in X$, $\cv p$ will just be $\set p$, and any singleton set is closed in a compact Hausdorff space. For any compact ordered space $(X,\le)$, since $\le$ is closed, $\cv p$ will also be closed for any $p\in X$. For complete lattices equipped with the canonical ultrastructures, we have:

\begin{lemma}\label{lem:cancvclosed}
  If $P$ is a complete lattice equipped with the canonical ultrastructure, then $\cv p \in \clc P$ for any $p\in P$.
\end{lemma}
\begin{proof}
  We verify directly. For any $f : S \to \cv p \subseteq P$ and any $\mu\in\beta S$,
  \[ \pair{f,\mu} = \bigvee_{A\in\mu} f^A. \]
  Since $f(S) \subseteq \cv p$, $f^A \le p$ for any $A\in\mu$. It follows that $\pair{f,\mu} \le p$, thus $\cv p$ is closed.
\end{proof}

On the other hand, $\dv p$ is always closed for any ultraposet:

\begin{lemma}
  Let $P$ be an ultraposet. Then for any $p\in P$, $\dv p \in \mathrm{Cl}_{\dv}(P)$.
\end{lemma}
\begin{proof}
  $\dv p$ is upward closed, thus we only need to show it is closed. Given any set $S$ and any $\mu\in\beta S$, consider the ultrapower $p^{\mu}$. We have shown that $p \le p^{\mu}$. Now for any family $f : S \to \dv p \subseteq P$, the constant function $p : S \to P$ sending everything to $p$ will be less than $f$ in the poset $[S,P]$, thus by functoriality of $\pair{-,\mu}$ we have
  \[ p \le p^{\mu} = \pair{p,\mu} \le \pair{f,\mu}. \]
  Hence, $\pair{f,\mu}\in\dv p$, and this shows $\dv p$ is closed.
\end{proof}

It is easy to see that, $\clc P$ and $\cld P$ again are closed under arbitrary intersection and finite union:

\begin{lemma}
  $\clc P$ and $\cld P$ are both closed under arbitrary intersection and finite union.
\end{lemma}
\begin{proof}
  Straight forward from Lemma~\Ref{lem:closedtop}, and the fact that both downward and upward closed subsets are closed under arbitrary intersection and union.
\end{proof}

Thus, $\clc P$ and $\cld P$ can also be viewed as topologies of closed sets on $P$, which are both coarser than $\cl P$. The nice thing now of course is that left and right ultrafunctors will be continuous for the topologies induced by $\clc-$ and $\cld-$, respectively:

\begin{lemma}\label{lem:lrultracont}
  A left (resp. right) ultrafunctor $\varphi : P \to Q$ is continuous when $P,Q$ are equipped with topologies $\clc P,\clc Q$ (resp. $\cld P,\cld Q$).
\end{lemma}
\begin{proof}
  Suppose $\varphi$ is a left ultrafunctor, and $L \in \clc Q$. Firstly, observe that $\varphi\inv(L)$ is indeed downward closed. We further show it is closed. Consider any $f : S \to \varphi\inv(L)$ and $\mu\in\beta S$. By assumption, $\varphi f : S \to L$, thus $\pair{\varphi f,\mu} \in L$. Since $\varphi$ is a left ultrafunctor,
  \[ \varphi\pair{f,\mu} \le \pair{\varphi f,\mu}. \]
  But $L$ is by definition downward closed, thus $\varphi\pair{f,\mu} \in L$, or equivalently, $\pair{f,\mu} \in \varphi\inv(L)$, hence $\varphi\inv(L)$ is also closed. The case for right ultrafunctors is completely symmetric.
\end{proof}

However, a left (resp. right) ultrafunctor viewed as a continuous function may not be closed anymore. This already fails on the level of preserving downward (resp. upward) closed subsets, since if $K$ is downward (resp. upward) closed in $P$, $\varphi(K)$ may \emph{not} be downward (resp. upward) closed. We then obtain two functors as follows,
\[ \clc- : \UltLPos \to \CTop, \quad \cld- : \UltRPos \to \CTop. \]

Furthermore, the construction $P \mapsto \clc P$ and $P \mapsto \cld P$, viewing $\clc P$ and $\cld P$ as \emph{lattices}, are both in some sense \emph{representable}:



\begin{theorem}\label{thm:lultclosedsubset}
  Left (resp. right) ultrafunctors from $P$ to $\mb 2$ contravariantly (resp. covariantly) correspond to closed and downward (resp. upward) closed subsets of $P$, i.e. there are natural isomorphisms of posets
  \[ \UltLPos(P,\mb 2)\op \cong \clc P, \quad \UltRPos(P,\mb 2) \cong \cld P. \]
\end{theorem}
\begin{proof}
  Given a left ultrafunctor $\varphi : P \to \mb 2$, since $\set 0$ is closed and downward closed in $\mb 2$, by Lemma~\ref{lem:lrultracont} $\varphi\inv(0) \in \clc P$. Given any closed and downward closed subsets $Q \subseteq P$, we show the induced map $\varphi : P \to \mb 2$ given by $\varphi(p) = 0$ iff $p\in Q$ is a left ultrafunctor. First of all, since $Q$ is downward closed, this is indeed monotone. To this end, consider $f : S \to P$ and $\mu\in\beta S$. We consider two cases:
  \begin{itemize}
  \item If $f\inv(Q) \in \mu$, since $Q$ is closed, $\pair{f,\mu} \in Q$. This means $\varphi\pair{f,\mu} = 0$, and hence we must have $\varphi\pair{f,\mu} \le \pair{\varphi f,\mu}$.
  \item If $f\inv(Q)\not\in\mu$, then $f\inv(P\backslash Q)\in\mu$. This suggests that the composite $\varphi f$ takes value 1 on the subset $f\inv(P\backslash Q)\in\mu$. Hence, by the ultrastructure on $\mb 2$, we have $\pair{\varphi f,\mu} = 1$, and we are done again.
  \end{itemize}
  Evidently the two constructions are inverses to each other. If $\varphi,\psi : P \to \mb 2$ are two left ultrafunctors with $\varphi \le \psi$, this implies $\psi\inv(0) \subseteq \varphi\inv(0)$, thus the isomorphism is contravariant.

  The case for right ultrafunctors is completely symmetric, taking a right ultrafunctor $\varphi : P \to \mb 2$ to the subset $\varphi\inv(1) \in \cld P$, and vice versa. Since $\varphi \le \psi$ implies $\varphi\inv(1) \subseteq \psi\inv(1)$, the isomorphism becomes covariant.
\end{proof}

Since an ultrafunctor is by definition exactly both a left and right ultrafunctor, from the above theorem it follows that ultrafunctors from $P$ to $\mb 2$ classify complemented pairs $(K,\ov K)$ of subsets in $P$, such that $K \in \clc P$ and $\ov K \in \cld P$; and $(K,\ov K) \le (L,\ov L)$ iff $K \supseteq L$. We use $\clcd P$ to denote this poset, and thus we have the following corollary:

\begin{corollary}\label{cor:ulclcd}
  For any ultraposet $P$, we have
  \[ \UltPos(P,\mb 2) \cong \clcd P. \]
\end{corollary}
\begin{proof}
  Straight froward from Theorem~\ref{thm:lultclosedsubset}. 
\end{proof}

\begin{remark}\label{rem:patchtop}
  The construction $\UltPos(-,\mb 2)$ in fact provides another functor from ultraposets to topological spaces:
  \[ \clcd- : \UltPos \to \CTop, \]
  where it sends each ultraposet $P$ to the socalled \emph{patch topology}, generated by a subbasis $\scomp{K,\ov K}{(K,\ov K) \in \clcd P}$. Since by definition for any such pair $(K,\ov K) \in \clcd P$, both $K$ and $\ov K$ lies in $\cl P$, hence is compact, the resulting space $P$ with the patch topology will again be compact.
\end{remark}

At this point, we have completed all the preliminaries, and are ready to discuss the embeddings from $\DL\op$ and $\CLoc$ to $\UltPos$ and $\UltLPos$, respectively. The following section first shows how to construct an ultraposet from a given distributive lattice or a coherent locale, by looking at the poset of \emph{models/points}. 

\section{Ultraposet from Distributive Lattice and Coherent Locale}\label{sec:uldiscoh}

The goal in this section is to construct two embedding functors,
\[ \Mod : \DL\op \to \UltPos, \quad \pt : \CLoc \to \UltLPos, \]
where $\Mod$ sends a distributive lattice to its poset of \emph{models}, and $\pt$ sends a coherent locale to its poset of \emph{points}. 

We start from a general \emph{poset}. For any poset $P$, by Example~\Ref{exm:catconsultrapos}, the poset $[P,\mb 2]$ of monotone maps from $P$ to $\mb 2$ will be an ultraposet. The description of the ultrastructure on $\mb 2$ also extends to these presheaf posets: For any $f : S \to [P,\mb 2]$ and any $\mu\in\beta S$, we have for $i \in \set{0,1}$,
\[ \pair{f,\mu}(p) = i \eff \scomp{s\in S}{f(s)(p) = i} \in \mu, \]

Notice that the poset $[P,\mb 2]$ is \emph{contravariantly isomorphic} to the poset of downward closed subsets of $P$, by identifying $x\in [P,\mb 2]$ with the preimage $x\inv(0)$ of 0. Henceforth, we will identify $[P,\mb 2]$ with the poset of downward closed subsets of $P$ under $\supseteq$. From this perspective, for any $p\in P$, we define a subset $C_p \subseteq [P,\mb 2]$ as follows,
\[ C_p := \scomp{I\in[P,\mb 2]}{p \in I}, \]
which is the set of all downward closed subsets that contain $p$. The ultrastructure on $[P,\mb 2]$ can now be described as follows,
\begin{align*}
  p \in \pair{f,\mu} \eff \scomp{s\in S}{p \in f(s)} \in \mu \eff f\inv(C_p) \in \mu
\end{align*}
Similarly, let $D_p$ be the \emph{complement} of $C_p$,
\[ D_p := \scomp{I}{p\not\in I}, \]
and we equivalently have
\[ \pair{f,\mu} = \scomp{p}{f\inv(D_p)\not\in\mu}. \]
Subsets of $[P,\mb 2]$ of the form $C_p,D_p$ will be called \emph{primitive}. They are particularly important because the pair $(C_p,D_p)$ lies in $\clcd{[P,\mb 2]}$:

\begin{lemma}\label{lem:primclosed}
  For any $p\in P$, $C_p \in \clc{[P,\mb 2]}$ and $D_p \in \cld{[P,\mb 2]}$. 
\end{lemma}
\begin{proof}
  Firstly, by definition $C_p$ is downward closed and $D_p$ is upward closed (recall that the order in $[P,\mb 2]$ is the \emph{converse} of inclusion of downward closed subsets in $P$). Suppose we have $f : S \to C_p \subseteq [P,\mb 2]$ and $\mu\in\beta S$, by the above characterisation,
  \[ \pair{f,\mu} = \scomp{q}{f\inv(C_q)\in\mu}. \]
  In particular, $f\inv(C_p) = S$, thus $p\in\pair{f,\mu}$. It follows that $\pair{f,\mu} \in C_p$, hence $C_p$ is also closed. The case for $D_P$ is similar.
\end{proof}

Given any monotone map $\sigma : P \to Q$, there is also an induced morphism from $[Q,\mb 2]$ to $[P,\mb 2]$, given by pre-composing with $\sigma$,
\[ \sigma^{*} : [Q,\mb 2] \to [P,\mb 2]. \]
Viewing an element $I$ in $[Q,\mb 2]$ as a downward closed subset of $Q$, we have
\[ \sigma^{*}(I) = \scomp{p}{\sigma p\in I} = \sigma\inv(I). \]
The inverse image of $\sigma^{*}$ actually preserves primitive subsets:

\begin{lemma}
\label{lem:prprimitive}
  For any monotone map $\sigma : P \to Q$, the inverse image $(\sigma^{*})\inv$ preserves primitive subsets. Concretly, for any $p\in P$,
  \[ (\sigma^{*})\inv(C_p) = C_{\sigma p}, \quad (\sigma^{*})\inv(D_p) = D_{\sigma p}. \]
\end{lemma}
\begin{proof}
  Given any $p \in P$ and $I \in [Q,\mb 2]$, by definition
  \[ I \in C_{\sigma p} \eff \sigma p \in I \eff p\in\sigma^{*}(I) \eff \sigma^{*}(I) \in C_p. \]
  This implies that 
  \[ (\sigma^{*})\inv(C_p) = C_{\sigma p}. \]
  Then we have
  \[ (\sigma^{*})\inv(D_p) = (\sigma^{*})\inv(\ov{C_p}) = \ov{(\sigma^{*})\inv(C_p)} = \ov{C_{\sigma p}} = D_{\sigma p}. \qedhere \]
\end{proof}

Lemma~\Ref{lem:prprimitive} actually implies that monotone functions between arbitrary posets induces ultrafunctors between their corresponding presheaf posets:
\begin{theorem}
\label{thm:posultra}
  The construction $P\mapsto[P,\mb 2]$ establishes a functor
  \[ [-,\mb 2] : \Pos\op \to \UltPos. \]
\end{theorem}
\begin{proof}
  We show for any monotone map $\sigma : P \to Q$, the induced morphism $\sigma^{*}:[Q,\mb 2] \to [P,\mb 2]$ is an ultrafunctor. Given any family $f : S \to [Q,\mb 2]$ and any $\mu\in S$, by definition
  \[ \pair{f,\mu} = \scomp{q}{f\inv(C_q) \in \mu}. \]
  Hence, we have
  \[ \sigma^{*}\pair{f,\mu} = \sigma\inv(\pair{f,\mu}) = \scomp{p}{f\inv(C_{\sigma p})\in\mu}. \]
  On the other hand, by Lemma~\Ref{lem:prprimitive},
  \begin{align*}
    \pair{\sigma^{*}f,\mu} &= \scomp{p}{(\sigma^{*}f)\inv(C_p)\in\mu} \\
                           &= \scomp{p}{f\inv(\sigma^{*})\inv(C_p) \in \mu} \\
                           &= \scomp{p}{f\inv(C_{\sigma p}) \in \mu}.
  \end{align*}
  Thus, $\sigma^{*}\pair{f,\mu}$ and $\pair{\sigma^{*}f,\mu}$ coincide, which implies $\sigma^{*}$ is indeed an ultrafunctor.
\end{proof}



Now we consider the case for distributive lattices. Given a distributive lattice $D$, we may view it as a \emph{propositional theory}. In this case, a \emph{model} $x$ of $D$ will simply be a distributive lattice morphism $x : D \to \mb 2$, viz. a compatible assignment of truth values of propositions in $D$. We use $\Mod(D)$ to denote the poset of models of $D$, under point-wise order. Notice that $\Mod(D)$ is by definition a subposet of $[D,\mb 2]$, and under our identification of $[D,\mb 2]$ with downward closed subset of $D$, $\Mod(D)$ is isomorphic to the set of \emph{prime ideals} of $D$, under $\supseteq$; see Lemma~\Ref{lem:modelprim} in Appendix~\Ref{app}.

To show $\Mod(D)$ is an ultraposet, we show it is \emph{closed} in $[D,\mb 2]$:

\begin{theorem}[\L os Ultraproduct Theorem for Propositional Logic]
\label{thm:los}
  $\Mod(D)$ is a closed subposet of $[D,\mb 2]$.
\end{theorem}
\begin{proof}
  Given a family of models $\set{x_s}_{s\in S}$ indexed by a set $S$, since the ultrastructure on $[D,\mb 2]$ is point-wise, their ultraproduct in $[D,\mb 2]$ is the following composite map,
  \[
    \xymatrix@C=3pc{ D \ar[r]^{\set{x_s}_{s\in S}} & \mb 2^S \ar[r]^{\pair{-,\mu}} & \mb 2 }
  \]
  $\mb 2^S$ has a point-wise lattice structure induced by that on $\mb 2$, thus the family $\set{x_s}_{s\in S} : D \to \mb 2^S$ will be a distributive lattice morphism, because each $x_s$ is. For the map $\pair{-,\mu}$, by the description of the ultrastructure on $\mb 2$ in Section~\Ref{sec:basic}, for any $x,y\in \mb 2^S$,
  \begin{align*}
    \pair{x\vee y,\mu} = 0
    &\eff (x\vee y)\inv(1) \not \in \mu \\
    &\eff x\inv(1) \cup y\inv(1) \not\in \mu \\
    &\eff x\inv(1) \not\in\mu \conjt y\inv(1)\not\in\mu \\
    &\eff \pair{x,\mu} = 0 \conjt \pair{y,\mu} = 0
  \end{align*}
  This implies that
  \[ \pair{x\vee y,\mu} = \pair{x,\mu} \vee \pair{y,\mu}. \]
  Similarly, we can show for conjunction that
  \[ \pair{x \wedge y,\mu} = 1 \eff \pair{x,\mu} = 1 \conjt \pair{y,\mu} = 1, \]
  Thus $\pair{-,\mu}$ preserves conjunction as well. It is also easy to see that $\pair{-,\mu}$ preserves the top and bottom element, thus it is a distributive lattice morphism. Hence, the composite will also be a model of $D$, and this shows that $\Mod(d)$ is closed under ultraproducts in $[D,\mb 2]$.
\end{proof}

Furthermore, the ultrafunctor between presheaf posets identified in Theorem~\Ref{thm:posultra} also respects models:

\begin{theorem}\label{thm:modelultra}
  The construction $D \mapsto \Mod(D)$ extends to a functor
  \[ \Mod : \DL\op \to \UltPos. \]
\end{theorem}
\begin{proof}
  We have seen from Theorem~\ref{thm:posultra} that the construction $[-,\mb 2]$ consists of a functor from $\Pos\op$ to $\UltPos$. By Theorem~\ref{thm:los}, for any $D\in\DL$, $\Mod(D) \subseteq [D,\mb 2]$ is closed, thus by Lemma~\Ref{lem:clores} we only need to show that for any distributive lattice homomorphism $\sigma : C \to D$, the induced morphism $\sigma^{*}$ takes models to models. But this is evident from definition, since $\sigma^{*}$ is given by pre-composition.
\end{proof}

Finally, we extend this picture to coherent locales. It is well-known that every coherent locale is isomorphic to one of the form $\idl(D)$ for some distributive lattice $D$, where $\idl(D)$ is the frame of ideals on $D$.

Furthermore, the poset of points of $\idl(D)$ will be isomorphic to $\Mod(D)$. A \emph{point} $x$ of $\idl(D)$ is equivalently a frame morphism
\[ x^{*} : \idl(D) \to \mb 2, \]
which preserves finite meets and arbitrary joins in $\idl(D)$. The distributive lattice $D$ also fully faithfully embeds into $\idl(D)$ as a sub-distributive lattice, given by
\[ p \mapsto \cv p. \]
Thus, by pre-composing with the embedding $D \hook \idl(D)$, we obtain a monotone function $\pt(\idl(D)) \to \Mod(D)$, and it is well-known that this is an \emph{isomorphism}. Its inverse takes a model $x$ of $D$, maps it to the frame morphism $x^{*} : \idl(D) \to \mb 2$ as follows,
\[ x^{*}(I) = \bigvee_{p\in I}x(p). \]
See~\cite[Ch. II.3]{johnstone1982stone} for a more detailed description on this.

Thus, for any coherent locale, we have already obtained an ultraposet by looking at its poset of points. We only need to show that this construction is functorial:

\begin{theorem}\label{thm:pointlultra}
  The assignment $L \mapsto \pt(L)$ extends to a functor
  \[ \pt(-) : \CLoc \to \UltLPos. \]
\end{theorem}
\begin{proof}
  We need to show for any localic map $g : \idl(C) \to \idl(D)$, the induced map $\pt(g) : \Mod(C) \to \Mod(D)$, which is given by pre-composing with the inverse image $g^{*}$ of $g$, is a left ultrafunctor. Consider any family $f : S \to \Mod(C)$ and $\mu\in\beta S$, for any $p\in D$,
  \[ (\pt(g)\pair{f,\mu})(p) = \pair{f,\mu}(g^{*}(\cv p)) = \bigvee_{q\in g^{*}\cv p}\pair{f,\mu}(q) = \bigvee_{q\in g^{*}(\cv p)}\bigvee_{A\in\mu}\bigwedge_{s\in A}fsq. \]
  On the other hand,
  \[ \pair{\pt(g)f,\mu}(p) = \bigvee_{A\in\mu}\bigwedge_{s\in A}(\pt(g)f)sq = \bigvee_{A\in\mu}\bigwedge_{s\in A}\bigvee_{q\in g^{*}(\cv p)}fsq. \]
  From the above description, it easily follows that
  \[ g\pair{f,\mu}(p) \le \pair{gf,\mu}(p),\ \forall p\in D, \]
  Hence, $\pt(g) : \Mod(C) \to \Mod(D)$ is a left ultrafunctor.
\end{proof}

\begin{remark}\label{rem:betaloc}
  Theorem~\Ref{thm:pointlultra} shows that the poset of points of each coherent locale is a closed subset of a presheaf ultraposet with a canonical ultrastructure. This actually holds for a wider class of locales. \cite{di2022geometry} has shown that this is more generally true for the class of locales which are right Kan injective w.r.t. embeddings of the form $X \hook \beta X$, where $X$ is viewed as a discrete space and $\beta X$ is its Stone-\u{C}ech compatification. 
\end{remark}

\section{Reconstruction from Models and Points}\label{sec:reconstruct}

In the previous section, we have successfully constructed two functors
\[ \Mod : \DL\op \to \UltPos, \quad \pt : \CLoc \to \UltLPos. \]
Our final goal is to show that they are \emph{fully faithful}. However, in this section, we first establish a partial result, by showing that we can reconstruct the coherent locale $\idl(D)$ from its ultraposet of points $\Mod(D)$. Concretely, we show there is a natural isomorphism of the following type,
\[ \UltLPos(\Mod(D),\mb 2) \cong \idl(D). \]
From this, we also prove that we can also reconstruct $D$ itself from $\Mod(D)$ by looking at ultrafunctors instead,
\[ \UltPos(\Mod(D),\mb 2) \cong D. \]

In Section~\Ref{sec:top}, we have already shown in Theorem~\Ref{thm:lultclosedsubset} that left ultrafunctors from any ultraposet $P$ to $\mb 2$ classify closed and downward closed subsets of $P$ in a contravariant way. Hence, we only need to show there is an isomorphism
\[ \idl(D)\op \cong \mathrm{Cl}_{\cv}(\Mod(D)). \]

We proceed by constructing a Galois connection between ideals on $D$ and subsets of $\Mod(D)$. Given an ideal $I$ in $D$, we consider the subset $K_I \subseteq \Mod(D)$ defined as follows,\footnote{Recall Section~\Ref{sec:uldiscoh} has identified $\Mod(D)$ with the set of prime ideals on $D$ under $\supseteq$.}
\[ K_I := \scomp{x\in\Mod(D)}{x \supseteq I}. \]

\begin{lemma}
  $K_I \subseteq \Mod(D)$ is closed and downward closed for any ideal $I$.
\end{lemma}
\begin{proof}
  $K_I$ is evidently downward closed. For a map $f : S \to K_I \subseteq \Mod(D)$ and any $\mu\in\beta S$, since $\Mod(D)$ is closed in $[D,\mb 2]$ by Theorem~\Ref{thm:los},
  \[ \pair{f,\mu} = \scomp{p}{f\inv(C_p) \in \mu} \in K_I \eff \scomp{p}{f\inv(C_p)\in\mu} \supseteq I. \]
  For any $s\in S$, since by assumption $fs\in K_I$, it follows that $fs \supseteq I$, thus for any $p\in I$ we have $fs \in C_p$. This means $f\inv(C_p) = S \in \mu$, and thus $\pair{f,\mu} \in K_I$. This proves $K_I \in \clc{\Mod(D)}$.
\end{proof}

To construct an ideal of $D$ from a subset of $\Mod(D)$, we first develop some preliminaries. For any $p\in D$, let $B_p$ and $O_p$ denote the intersection $C_p \cap \Mod(D)$ and $D_p \cap \Mod(D)$, respectively. Since $\Mod(D)$ is a closed subset of $[D,\mb 2]$, and by Lemma~\Ref{lem:primclosed} $(C_p,D_p) \in \clcd{[D,\mb 2]}$, it follows that $(B_p,O_p) \in \clcd{\Mod(D)}$. Furthermore, $B_p$ and $O_p$ can also be used to describe the ultrastructure on $\Mod(D)$: For any family $f : S \to \Mod(D)$ and any $\mu\in\beta S$,
\[ \pair{f,\mu} = \scomp{p}{f\inv(C_p)\in\mu} = \scomp{p}{f\inv(B_p)\in\mu} = \scomp{p}{f\inv(O_p)\not\in\mu}. \]

These subsets $B_p$ and $O_p$ also interact well with the lattice structure:

\begin{lemma}
\label{lem:primitivestr}
  For any $p,q \in D$, we have $B_{p\vee q} = B_p \cap B_q$ and $B_{p\wedge q} = B_p \cup B_q$. Dually, $O_{p\vee q} = O_p \cup O_q$ and $O_{p\wedge q} = O_p \cap O_q$.
\end{lemma}
\begin{proof}
  We prove the case for $p\vee q$. By definition,
  \[ B_{p\vee q} = \scomp{x}{x(p\vee q) = 0} = \scomp{x}{x(p) = 0 \conjt x(q) = 0} = B_p \cap B_q. \]
  Similarly for conjunction, and the case for $O_p$ is completely symmetric.
\end{proof}

This way, given any subset $P \subseteq \Mod(D)$ of $\Mod(D)$, we can construct an ideal $I_P$ of $D$ as follows,
\[ I_P := \scomp{p}{P \subseteq B_p}. \]

\begin{corollary}
  $I_P$ is an ideal of $D$ for any subset $P \subseteq \Mod(D)$.
\end{corollary}
\begin{proof}
  This is essentially due to Lemma~\ref{lem:primitivestr}. By definition, for any $p\le q$, $B_q \subseteq B_p$. Hence, if $q\in I_P$ then $P \subseteq B_q \subseteq B_p$, which implies $p \in I_P$ as well, hence $I_P$ is downward closed. For any $p,q\in I_P$, $P \subseteq B_p$ and $P \subseteq B_q$ implies $P \subseteq B_p \cap B_q = B_{p \vee q}$, hence $p \vee q \in I_P$ as well. This proves $I_P$ is an ideal.
\end{proof}

Our task is then to show that the constructions $I \mapsto K_I$ and $P \mapsto I_P$ consist of a contravariant isomorphism between $\idl(D)$ and $\clc{\Mod(D)}$. One direction is actually quite straight forward:

\begin{lemma}
  The construction $I \mapsto K_I \mapsto I_{K_I}$ returns $I$ itself.
\end{lemma}
\begin{proof}
  By definition, for any ideal $I$ in $D$,
  \begin{align*}
    I_{K_I} &= \scomp{p}{K_I \subseteq B_p} = \scomp{p}{\forall x\in K_I,\,x\in B_p} \\
            &= \scomp{p}{\forall x\in\Mod(D),\,x \supseteq I \nt p\in x} \\
            &= \bigcap \scomp{x\in\Mod(D)}{x \supseteq I} = I.
  \end{align*}
  The final equality is due to the well-known fact that any ideal in a distributive lattice is equal to the intersection of all prime ideals extending it. For a proof, see Appendix~\Ref{app}, Lemma~\ref{lem:idealintersecprime}. 
\end{proof}

To show for any $K \in \clc{\Mod(D)}$, $K_{I_K}$ also gives back $K$ itself, we need some further observations.

\begin{lemma}\label{lem:primdis}
  For any poset $P$ and any $I \in [P,\mb 2]$, we have\footnote{Notice that we have identified $[P,\mb 2]$ as the set of downward closed subsets of $P$ under $\supseteq$, thus $\cv I$ actually represents the set of all downward closed subsets that \emph{contains} $I$, and similarly $\dv I$ represents the set of all downward closed subsets that is contained in $I$.}
  \[ \cv I = \bigcap_{p\in I}C_p, \quad \dv I = \bigcap_{p\not\in I}D_p. \]
\end{lemma}
\begin{proof}
  First of all, if $J \supseteq I$, then for any $p\in I$, we also have $p\in J$, thus $J \in \bigcap_{p\in I}C_p$. On the other hand, for any $J$ that $J \supsetneq I$, there exists some $p\in I$ that $p\not\in J$. It follows that $J \not\in C_p$, thus $\cv I = \bigcap_{p\in I}C_p$. The proof $\dv I = \bigcap_{p\not\in I}D_p$ is completely similar.
\end{proof}

Lemma~\Ref{lem:primdis} implies that primitive subsets of $[P,\mb 2]$ are sufficient to distinguish all elements in $[P,\mb 2]$. If moreover we have a distributive lattice $D$, subsets of the form $B_p$ and $O_p$ also suffices:

\begin{lemma}\label{lem:disprimdis}
  Let $D$ be a distributive lattice, and let $K \in \clc{\Mod(D)}$. If $x\in\Mod(D)$ does \emph{not} lie in $K$, then there exists some $p\in D$ such that $K \subseteq B_p$ and $x \in O_p$.
\end{lemma}
\begin{proof}
  Since $K$ is downward closed and $x \not\in K$, it follows that $\dv x$ is disjoint from $K$. By Lemma~\ref{lem:primdis}, it follows that $K \cap \bigcap_{p\not\in x}O_p = \emptyset$. Since both $K$ and all $O_p$ are closed, and $\Mod(D)$ is compact by Lemma~\ref{lem:compact}, there must exists some finite subfamily that
  \[ K \cap \bigcap_{i=1}^nO_{p_i} = \emptyset. \]
  Notice that $K$ must be in this finite subfamily, because by definition $x\in\bigcap_iO_{p_i}$. Then let $p_1\wedge\cdots\wedge p_n$ be $p$. By Lemma~\Ref{lem:primitivestr},
  \[ x \in O_{p_1} \cap \cdots \cap O_{p_n} = O_{p_1\wedge\cdots\wedge p_n} = O_p, \]
  and $K$ is disjiont from $O_p$. Since $B_p$ and $O_p$ are complemented, $K \subseteq B_p$.
\end{proof}

\begin{corollary}\label{cor:clcintersec}
  For any $K\in \mathrm{Cl}_{\cv}(\Mod(D))$, we have $K = \bigcap_{K\subseteq B_p}B_p$.
\end{corollary}
\begin{proof}
  For any $x\not\in K$, by the previous lemma there exists some $p\in D$ that $K \subseteq B_p$ and $x\not\in B_p$. Hence, $K = \bigcap_{K\subseteq B_p}B_p$.
\end{proof}

\begin{corollary}
  For any $K\in \mathrm{Cl}_{\cv}(\Mod(D))$, the construction $K \mapsto I_K \mapsto K_{I_K}$ gives back to itself.
\end{corollary}
\begin{proof}
  By definition, we have
  \begin{align*}
    K_{I_K} &= \scomp{x\in\Mod(D)}{x \supseteq I_K} \\
            &= \scomp{x\in\Mod(D)}{\forall p\in D,\,K\subseteq B_p \nt x\in B_p} \\
            &= \bigcap_{K\subseteq B_p}B_p = K. \qedhere
  \end{align*}
\end{proof}

Hence, indeed we have a bijection bewteen $\idl(D)$ and $\clc{\Mod(D)}$. It is evident from our definition of $I_{(-)}$ and $K_{(-)}$ that this bijection is contravariant, thus we have the following result:

\begin{theorem}\label{thm:recidl}
  For any distributive lattice $D$, we have the following lattice isomorphisms,
  \[ \UltLPos(\Mod(D),\mb 2) \cong \mathrm{Cl}_{\cv}(\Mod(D))\op \cong \idl(D). \]
\end{theorem}

\begin{remark}
  We also have a completely dual statement of Theorem~\Ref{thm:recidl}. Let $\fil(D)$ denote the lattice of \emph{filters} on $D$. By essentially the same argument as before, we can show the following isomorphisms,
  \[ \UltRPos(\Mod(D),\mb 2) \cong \cld{\Mod(D)} \cong \fil(D). \]
  This partially reflects the advantage of the approach of ultraposets, since it encompasses all there results with a unifying theoretic framework.
\end{remark}

Now it is straight forward to derive the reconstruction for the distributive lattice itself, viz. to show $D$ is isomorphic to ultrafunctors from $\Mod(D)$ to $\mb 2$. Under the embedding $D \hook \idl(D)$, with each $p\in D$ mapping to $\cv p$, the isomorphism $\idl(D) \cong \UltLPos(\Mod(D),\mb 2)$ by definition takes it to the following closed subset $K_{\cv p}$ in $\Mod(D)$:
\[ K_{\cv p} = \scomp{x\in\Mod(D)}{x \supseteq \cv p} = \scomp{x\in\Mod(D)}{p \in x} = B_p. \]
As commented before, $(B_p,O_p) \in \clcd{\Mod(D)}$, and we have seen from Corollary~\Ref{cor:ulclcd} that such pairs exactly correspond to ultrafunctors from $\Mod(D)$ to $\mb 2$, the isomorphism $\idl(D) \cong \UltLPos(\Mod(D),\mb 2)$ indeed restricts to an embedding
\[ D \hook \UltPos(\Mod(D),\mb 2). \]
To show this is an isomorphism, we only need to prove that every pair $(K,\ov K)\in\clcd{\Mod(D)}$ comes from some $p\in D$.

It is also well-known that the embedding $D \hook \idl(D)$ identifies $D$ as the sublattice of \emph{finite element} in $\idl(D)$ (cf.~\cite[Ch. II.3]{johnstone1982stone}). An element $a$ in a complete lattice $P$ is \emph{finite} if for any ideal $I$ of $P$, $a \le \bigvee I$ iff $a \in I$. In the lattice $\clc{\Mod(D)}\op$, an element $K$ is finite iff for any ideal $\set{K_i}_{i\in I}$ in $\clc{\Mod(D)}\op$, if $K = \bigcap_{i\in I}K_i$, then $K = K_i$ for some $i\in I$.

\begin{lemma}\label{lem:clcomclfinite}
  For any $K\in\clc{\Mod(D)}$, if its complement is also closed, viz. $\ov K \in \cld{\Mod(D)}$, then it is finite.
\end{lemma}
\begin{proof}
  Let $K\in\clc{\Mod(D)}$ with $\ov K \in \cld{\Mod(D)}$. Suppose $K = \bigcap_{i\in I}K_i$, with $\set{K_i}_{i\in I}$ an ideal of $\clc{\Mod(D)}\op$. If $K \subsetneq K_i$ for all $i\in I$, then $\ov K \cap K_i \neq\emptyset$ for all $i\in I$. It follows that for any finite subfamily of $\set{\ov K} \cup \set{K_i}_{i\in I}$, it has non-empty intersection, because $\set{K_i}_{i\in I}$ is an ideal. By compactness from Lemma~\Ref{lem:compact}, $\ov K \cap \bigcap_{i\in I}K_i$ is non-empty, contradictory to the fact that $K = \bigcap_{i\in I}K_i$. Thus, we must have $\ov K \cap K_i = \emptyset$, which implies $K = K_i$, for some $i\in I$. Hence $K$ is finite. 
\end{proof}

\begin{corollary}\label{cor:reul}
  The isomorphism $\UltLPos(\Mod(D),\mb 2) \cong \idl(D)$ restricts to an isomorphism
  \[ \UltPos(\Mod(D),\mb 2) \cong \clcd{\Mod(D)} \cong D. \]
\end{corollary}
\begin{proof}
  By Lemma~\Ref{lem:clcomclfinite}, every ultrafunctor from $\Mod(D)$ to $\mb 2$ classfies some finite element in $\clc{\Mod(D)}\op$, thus corresponds to some $p\in D$. 
\end{proof}

\section{Duality for Distributive Lattices and Coherent Locales}\label{sec:dual}

In this section, we furthermore show that the functors $\Mod$ and $\pt$ from $\DL\op$ and $\CLoc$ to $\UltPos$ and $\UltLPos$ are both fully faithful. We start by observing that the assignment $P \mapsto \UltLPos(P,\mb 2) \cong \clc P$ actually extends to a functor
\[ \Omega : \UltLPos \to \Loc, \]
which will in some sense be a \emph{left adjoint} of $\pt$. Then the reconstruction result of Theorem~\ref{thm:recidl} implies this adjunction restricts to an equivalence when restricted to $\CLoc$ and a suitable full subcategory of $\UltLPos$.

We've already seen that $\UltLPos(P,\mb 2)$ for any ultraposet $P$ is isomorphic to $\clc{\Mod(D)}\op$. From Lemma~\Ref{lem:closedtop}, we know $\clc{\Mod(D)}$ is closed under arbitrary intersection and finite union, thus $\clc{\Mod(D)}\op$ will be a locale. It is easy to see that this extends to a functor:

\begin{lemma}\label{lem:lulposetloc}
  The assignment $P \mapsto \UltLPos(P,\mb 2)$ exhibits a functor
  \[ \Omega : \UltLPos \to \Loc. \]
\end{lemma}
\begin{proof}
  We only need to verify that for any left ultrafunctor $\varphi : Q \to P$, precomposition with $\varphi$ induces a frame homomorphism as follows,
  \[ \varphi^{*} : \Omega(Q) \to \Omega(P). \]
  But this is evident, because from Lemma~\Ref{lem:ultracont} we already know that $\varphi^{*}$ simply takes $K \in \clc Q$ to $\varphi\inv(K) \in \clc P$, and $\varphi\inv$ preserves arbitrary intersection and union. 
\end{proof}

Let $\CUltLPos$ be the full subcategory of $\UltLPos$ which lies in the essential image of the functor $\pt : \CLoc \to \UltLPos$. Firstly, from Theorem~\Ref{thm:recidl}, it follows that $\Omega$ restricts to a functor of the form
\[ \Omega : \CUltLPos \to \CLoc. \]
Verifying $\Omega$ and $\pt$ forms an adjunction is straight forward. We can in fact show something more general:

\begin{lemma}\label{lem:ptomegaadj}
  For any ultraposet $P$ and coherent locale $L$, there is a natural isomorphism
  \[ \UltLPos(P,\pt(L)) \cong \Loc(\Omega(P),L). \]
\end{lemma}
\begin{proof}
  Suppose $L \cong \idl(D)$, consider the following commutative diagramme,
  \[
    \xymatrix{
      \UltLPos(P,\pt(L)) \ar@{^{(}->}[d] \ar[r]^<<<<{\Omega} & \Loc(\Omega(P),\Omega(\pt(L)) \cong L) \ar@{^{(}->}[d] \\
      \UltLPos(P,[D,\mb 2]) \ar[r]_{\cong} & [D,\Omega(P)]
    }
  \]
  The top arrow uses Lemma~\Ref{thm:recidl}. The left arrow is given by the embedding
  \[ \pt(L) \cong \Mod(D) \hook [D,\mb 2]. \]
  The right arrow sends $g : \Omega(P) \to L$ in $\Loc$ to the map $p \mapsto g^{*}(\cv p)$ from $D$ to $\Omega(P)$. The bottom arrow is in fact an isomorphism by definition, because $[D,\mb 2]$ is equipped with the point-wise ultrastructure. It is straight forward to verify the above diagramme is a pullback.
\end{proof}

\begin{remark}\label{rem:interconnecttopulloc}
  Lemma~\ref{lem:ptomegaadj} suggests that there should be an adjunction bewteen the category $\UltLPos$ of ultraposets with left ultrafunctors, and some full subcategory of locales, with $\Omega \dashv \pt$. This reminds us with the adjunction between locales and topological spaces, which suggestively is also written as $\Omega \dashv \pt$, as follows,
  \[ \Omega : \Top \gb \Loc : \pt. \]
  It is interesting to characterise what is exactly the subcategory of $\Loc$ that makes the adjunction works. We will say a little bit more about this in Section~\Ref{sec:concfuture}. 
\end{remark}

The following theorem is an immediate consequence, which is one of our main goals of this paper:

\begin{theorem}\label{thm:cloccluteq}
  $\Omega$ and $\pt$ restricts to an equivalence $\CLoc \simeq \CUltLPos$.
\end{theorem}
\begin{proof}
  By definition of $\CUltPos$ and Theorem~\Ref{thm:recidl}, the result in Lemma~\Ref{lem:ptomegaadj} implies that there is an adjunction as follows,
  \[ \Omega : \CUltLPos \gb \CLoc : \pt. \]
  Theorem~\Ref{thm:recidl} also implies that the counit $\Omega\circ\pt \nt \id$ of the adjunction is an isomorphism, hence $\pt$ must be fully faithful. By definition of $\CUltLPos$, $\pt$ is also essentially surjective, thus they form an equivalence.
\end{proof}

As usual, our next step is to restrict the above equivalence to an equivalence on distributive lattices. Just like $\UltLPos(-,\mb 2)$ exhibits a functor from $\UltLPos$ to $\Loc$, $\UltPos(-,\mb 2)$ similarly will establish a functor from $\UltPos$ to $\DL\op$. Let $\CUltPos$ be the full subcategory of $\UltPos$ that lies in the image of $\Mod : \DL\op \to \UltPos$. $\CUltPos$ has the same objects as $\CUltLPos$, only the morphisms are restricted to ultrafunctors. We then have:

\begin{lemma}\label{lem:omegarestrict}
  The functor $\UltPos(-,\mb 2)$ defines a functor from $\UltPos$ to $\DL\op$, which we denote as $\Omega_u$,
  \[ \Omega_u : \UltPos \to \DL\op. \]
  Furthermore, the following diagramme commute,
  \[
    \xymatrix{
      \CUltPos \ar@{>->}[d] \ar[r]^{\Omega_u} & \DL\op \ar@{>->}[d] \\
      \CUltLPos \ar[r]_{\Omega} & \Loc
    }
  \]
\end{lemma}
\begin{proof}
  Under the isomorphism given by Corollary~\Ref{cor:ulclcd}, we have
  \[ \UltPos(P,\mb 2) \cong \clcd P. \]
  Since both $\clc P$ and $\cld P$ are closed under finite intersection and union, it follows that $\clcd P$ will be a distributive lattice. Furthemore, under this isomorphism, any ultrafunctor $\varphi : P \to Q$ is mapped to the inverse image $\varphi\inv$, thus preserving the distributive lattice structure.
  
  Now for any ultraposet of the form $\Mod(D)$ for some distributive lattice $D$, we have already shown in Theorem~\Ref{thm:recidl} and Corollary~\Ref{cor:reul} that
  \[ \UltPos(\Mod(D),\mb 2) \cong D, \quad \UltLPos(\Mod(D),\mb 2) \cong \idl(D). \]
  It follows the above diagramme commutes on objects. It commutes on morphisms as well, because under the isomorphisms given by Theorem~\Ref{thm:lultclosedsubset} and Corollary~\Ref{cor:ulclcd}, both are mapped to the inverse image $\varphi\inv$. 
\end{proof}

Lemma~\Ref{lem:omegarestrict} now suggests that the two functors only differ by their domains. It immediately completes our task in this section:

\begin{theorem}\label{thm:disculeq}
  $\Mod : \DL\op \to \CUltPos$ is an equivalence.
\end{theorem}
\begin{proof}
  For any distributive lattices $C,D$, we have the following commuting diagramme by Lemma~\Ref{lem:omegarestrict} and Theorem~\Ref{thm:cloccluteq},
  \[
    \xymatrix{
      \UltPos(\Mod(C),\Mod(D)) \ar@<.5ex>[r]^<<<<<<{\Omega_u = \Omega} \ar@{^{(}->}[d] & \DL(D,C) \ar@{^{(}->}[d] \ar@<.5ex>[l]^>>>>>>{\Mod} \\
      \UltLPos(\Mod(C),\Mod(D)) \ar[r]_<<<<{\cong} & \Loc(\idl(C),\idl(D))
    }
  \]
  It follows that $\Mod$ is fully faithful, and essentially surjective by definition of $\CUltPos$, thus is an equivalence.
\end{proof}

In this section, we have mainly shown that the following two embedding functors are indeed \emph{fully faithful},
\[ \Mod : \DL\op \to \UltPos, \quad \pt : \CLoc \to \UltLPos, \]
and we have defined the essential images of these embeddings as $\CUltPos$ and $\CUltLPos$, respectively. However, we would also like an intrinsic characterisation of those ultraposets that lie in these essential images. This is the task of the next section.

\section{Separation Properties for Ultraposets}\label{sec:sep}

We start by discussing a bit more about the relation between ultraposet and distributive lattices in general. From Lemma~\Ref{lem:omegarestrict}, for any ultraposet $P$, $\Omega_u(P) \cong \UltPos(P,\mb 2)$ will be a distributive lattice. There is always an evaluation map $\eta : P \to \Mod(\Omega_u(P))$, given by
\[ \eta(p) = \varphi \mapsto \varphi(p). \]
This evaluation map is always an ultrafunctor:

\begin{lemma}\label{lem:evultra}
  For any ultraposet $P$, the evaluation map is an ultrafunctor,
  \[ \eta : P \to \Mod(\Omega_u(P)). \]
\end{lemma}
\begin{proof}
  For any family $f : S \to P$, $\mu\in\beta S$, and $\varphi \in \Omega_u(P)$, by definition
  \[ (\eta\pair{f,\mu})(\varphi) = \varphi\pair{f,\mu} = \pair{\varphi f,\mu}. \]
  Recall that the ultrastructure on $\Mod(\Omega_u(P))$ is inherited from the canonical ultrastructure on $[\Omega_u(P),\mb 2]$, which is point-wise, thus
  \[ \pair{\eta \circ f,\mu}(\varphi) = \pair{(\eta \circ f)(\varphi),\mu} = \pair{\varphi f,\mu}. \]
  This implies $\eta$ is an ultrafunctor.
\end{proof}

\begin{theorem}\label{thm:uldladj}
  $\Omega_u \dashv \Mod$ forms an adjunction as follows,
  \[ \Omega_u : \UltPos \gb \DL\op : \Mod, \]
  which identifies $\DL\op$ as a reflexive subcategory of $\UltPos$.
\end{theorem}
\begin{proof}
  Lemma~\Ref{lem:evultra} provides the unit, and Corollary~\Ref{cor:reul} gives the counit. Verifying they satisfy all the triangular identities is straight forward, thus $\Omega_u \dashv \Mod$ forms an adjunction. The counit by Corollary~\Ref{cor:reul} is a natural isomorphism, thus $\DL\op$ is a reflexive subcategory of $\UltPos$.
\end{proof}

From Theorem~\Ref{thm:uldladj}, it follows that the category $\CUltPos$ consists of exactly those ultraposets where the unit $\eta$ is an \emph{isomorphism}. Hence, our goal is then to characterise those ultraposets $P$ such that $\eta : P \cong \Mod(\Omega_u(P))$.

We start by observing some important properties of ultraposets of the form $\Mod(D)$ for some distributive lattice $D$. Recall that one important fact that has been used extensively in the last two sections is Lemma~\Ref{cor:clcintersec}, i.e. any $K \in \clc{\Mod(D)}$ can be written as an intersection of primitive subsets $K = \bigcap B_p$, where $(B_p,O_p)$ lies in $\clcd{\Mod(D)}$. This leads us to the following definition:

\begin{definition}\label{def:zerodim}
  Let $P$ be an ultraposet. We say $P$ is \emph{zero-dimensional} if for any $p,q\in P$, if $p \not\le q$, then there exists $(K,\ov K) \in \clcd{P}$, that $p \in \ov K$ and $q \in K$. 
\end{definition}

\begin{remark}\label{rem:zerobool}
  For an \emph{ultraset} $X$, an ultrafunctor from $X$ to $\mb 2$ is equivalently a pair of clopen subsets of $X$. Then, our notion of being zero-dimensional as an ultraposet reduces to the topological notion of a zero-dimensional space, or a Stone space, under the isomorphism $\UltSet \cong \Comp$. Thus, once we have shown $\CUltPos$ is exactly the category of zero-dimensional ultraposets, the Stone duality $\Bool\op \simeq \Stone$ will follow as a consequence, because the order on $\Mod(D)$ is discrete iff $D$ is Boolean. 
\end{remark}

One immediate observation is that any ultraposet of the form $\Mod(D)$ for some distributive lattice $D$ is zero-dimensional:

\begin{lemma}\label{lem:zeroinclcul}
  Any ultraposet in $\CUltPos$ is zero-dimensional.
\end{lemma}
\begin{proof}
  For any $x\in\Mod(D)$, since $\Mod(D)$ is a closed subset of $[D,\mb 2]$, and by Lemma~\ref{lem:cancvclosed} $\cv x$ is closed in $[D,\mb 2]$, it follows that $\cv x \in \clc{\Mod(D)}$ as well. If $x \not\le y$, then $y\not\in\cv x$, and it follows from Lemma~\Ref{cor:clcintersec} that there must exists some $B_P$ satisfying $\cv x \subseteq B_p$ and $y\not\in B_p$, viz. $y\in O_p$. Hence, $\Mod(D)$ is zero-dimensional.
\end{proof}

To show the other direction, as mentioned previously, we only need to show that the unit $\eta : P \to \Mod(\Omega_u(P))$ is in fact an isomorphism for any zero-dimensional ultraposet. The key observation is that any zero-dimensional ultraposet $P$ can in fact be viewed as a \emph{compact ordered space}, i.e. it can in fact be viewed as a compact Hausdorff space with a closed partial order, where its ultrastructure is provided by its topology:

\begin{lemma}\label{lem:zerotop}
  Let $P$ be an zero-dimensional ultraposet $P$. For any $f : S \to P$, $\mu : T \to \beta S$ and any $\nu\in\beta T$,
  \[ \pair{f,\pair{\mu,\nu}} = \pair{\pair{f,\mu},\nu}. \]
\end{lemma}
\begin{proof}
  Firstly, we show that if $(K,\ov K) \in \clcd P$, then for any $f : S \to P$ and any $\mu\in\beta S$, we must have
  \[ \pair{f,\mu} \in K \eff f\inv(K) \in \mu. \]
  The right to left is evident, because $K$ is closed in $P$. For the left to right direction, we use the fact that $\ov K$ is also closed, thus
  \[ f\inv(K)\not\in\mu \nt f\inv(\ov K) \in \mu \nt \pair{f,\mu}\in\ov K \nt \pair{f,\mu} \not\in K. \]
  Now suppose $P$ is zero-dimensional, and suppose for some $f,\mu,\nu$, we have
  \[ \pair{f,\pair{\mu,\nu}} < \pair{\pair{f,\mu},\nu}, \]
  where the inequality holds strictly. Hence, we can find $(K,\ov K)\in\clcd P$ that
  \[ \pair{f,\pair{\mu,\nu}} \in K, \quad \pair{\pair{f,\mu},\nu} \in \ov K. \]
  Under the previous equivalence, from $\pair{f,\pair{\mu,\nu}}\in K$ we know that
  \[ f\inv(K) \in \pair{\mu,\nu} \eff \scomp{t}{f\inv(K)\in\mu_t} \in \nu. \]
  From $\pair{\pair{f,\mu},\nu} \in \ov K$ we also know that
  \[ \pair{f,\mu}\inv(\ov K) \in \nu \eff \scomp{t}{\pair{f,\mu_t}\in\ov K} \in \nu \eff \scomp{t}{f\inv(\ov K) \in \mu_t} \in\nu. \]
  However, since $K$ and $\ov K$ are complemented, it follows that for any $t\in T$,
  \[ f\inv(K) \in \mu_t \eff f\inv(\ov K) \not\in \mu_t. \]
  Then $\nu$ contains two complemented subset of $T$, which cannot happen. Thus, we must have
  \[ \pair{f,\pair{\mu,\nu}} = \pair{\pair{f,\mu},\nu}. \qedhere \]
\end{proof}

Lemma~\Ref{lem:zerotop} implies that even if we consider the discrete order $P^d$ on $P$, the ultrastructure on $P$ will still be an ultrastructure on $P^d$. Under the isomorphism $\UltSet \cong \Comp$, it follows that the topology $\cl P$ of all closed sets on $P$ is \emph{compact Hausdorff}, and the ultrastructure on $P$ is identical to the ultrastructure provided by this compact Hausdorff topology.

Recall that a compact Hausdorff space equipped with an order is precisely a \emph{Priestley space}:

\begin{definition}
  A \emph{Priestley space} is a compact Hausdorff topological space $X$ equipped with a partial order $\le$, such that for any $x,y\in X$, if $x \not\le y$, then there exists a clopen downset $K$ in $X$ that $y \in K$ and $x \not\in K$. A morphism between Priestley spaces is simply a continuous map that is also monotone. 
\end{definition}

In a Priestley space, the order is necessarily closed. We will use $\Pries$ to denote the (2-)category of Pristley spaces, where the order between arrows is given point-wise. It is now immediate that there is an isomorphism of categories as follows: 

\begin{corollary}\label{cor:zeropriestley}
  The (2-)categories of Priestley spaces and zero-dimensional ultraposets with ultrafunctors between them are isomorphic.
\end{corollary}
\begin{proof}
  By Lemma~\Ref{lem:compord} and Lemma~\Ref{lem:zerotop}, a Priestley space $(X,\le)$ is exactly an zero-dimensional ultraposet. Given such $P,Q$, a function $\varphi : P \to Q$ is an ultrafunctor iff it is monotone, and preserves the ultrastructure. Notice that $\varphi$ preserves the ultrastructure, under the isomorphism $\UltSet \cong \Comp$, iff it is continuous for the topology $\cl P$ and $\cl Q$. This implies $\varphi$ is an ultrafunctor iff it is a morphism between Priestley spaces.
\end{proof}

\begin{remark}
  There are two well-known categories of certain topological structure isomorphic to $\Pries$, that of \emph{spectral spaces} and \emph{pairwise Stone spaces} (cf.~\cite{bezhanishvili2010bitopological,cornish1975h}). Given a Priestley space $(X,\le)$, there are two further topologies $\tau^u,\tau^d$ associated to it, where they are given by the topology of open upsets and open downsets, respectively. The set $X$ equipped with the topology $\tau^u$ makes it into a spectral space, and $X$ equipped with $\tau^u,\tau^d$ makes it into a pairwise Stone space. Viewing $X$ as a zero-dimensional ultraposet, the closed sets of these additional topologies are exactly given by $\clc X$ and $\cld X$, which naturally arise from the functor $\UltLPos(-,\mb 2)$ and $\UltRPos(-,\mb 2)$ by Theorem~\Ref{thm:lultclosedsubset}. The theory of ultraposets and (left/right) ultrafunctors again provides a unifying framework for studying these different topological notions. 
\end{remark}

The fact that the category $\Pries$ of Priestley spaces is dual to the category $\DL$ of distributive lattices is a classical result in duality theory; for a textbook account, see~\cite{gehrke2022topological}. In particular, every Priestley space is isomorphic to the space of models of its clopen upsets, thus we have:

\begin{corollary}\label{cor:zerounitiso}
  The unit $\eta : P \to \Mod(\Omega_u(P))$ is an isomorphism for all zero-dimensional ultraposet.
\end{corollary}
\begin{proof}
  Follows from the fact that $P$ is a Priestley space.
\end{proof}

\begin{theorem}\label{thm:culdl}
  $\CUltPos$ (resp. $\CUltLPos$) is exactly the category of zero-dimensional ultraposets with ultrafunctors (resp. left ultrafunctors) between them.
\end{theorem}
\begin{proof}
  The two categories $\CUltPos$ and $\CUltLPos$ have the same objects by construction. The remaining follows directly follows from Lemma~\Ref{lem:zeroinclcul} and Corollary~\Ref{cor:zerounitiso}. 
\end{proof}

As we've mentioned in Remark~\ref{rem:zerobool}, the notion of zero-dimensional ultraposet, when applied to ultrasets, gives exactly Stone spaces. Thus, the classical Stone duality follows as an easy consequence:

\begin{corollary}
  There is an equivalence of categories $\Stone \simeq \Bool\op$.
\end{corollary}
\begin{proof}
  For any distributive lattice $D$, $\Mod(D)$ is discrete iff $D$ is Boolean (cf.~\cite[Ch. II, 4.4]{johnstone1982stone}). Thus, by Theorem~\Ref{thm:disculeq} and Theorem~\Ref{thm:culdl}, $\Bool\op$ is equivalent to discrete zero-dimensional ultraposets, viz. zero-dimensional ultrasets, which, from Remark~\Ref{rem:zerobool}, are exactly Stone spaces.
\end{proof}

\section{Conclusion and Future Work}\label{sec:concfuture}

In this paper, we have provided a complete description of the following diagramme,
\[
  \xymatrix{
    \Bool\op \ar@{^{(}->}[d] \ar[r]^{\Mod}_{\simeq} & \Stone \ar@{^{(}->}[d] \ar@{^{(}->}[r] & \UltSet \ar@{^{(}->}[d] \\
    \DL\op \ar@{>->}[d] \ar[r]^{\Mod}_{\simeq} & \CUltPos \ar@{>->}[d] \ar@{^{(}->}[r] & \UltPos \ar@{>->}[d] \\
    \CLoc \ar[r]^{\pt}_{\simeq} & \CUltLPos \ar@{^{(}->}[r] & \UltLPos
  }
\]
In particular, following Lurie~\cite{lurie2018ultracategories}, we have defined (2-)categories $\UltPos$, $\UltLPos$ and $\UltRPos$, of ultraposets with (left/right) ultrafunctors, and showed that the opposite of the category of distributive lattices $\DL\op$ fully faithfully embeds into $\UltPos$, and the category of coherent locales $\CLoc$ fully faithfully embeds into $\UltLPos$. The essential image of these embeddings corresponds to what we call zero-dimensional ultraposets.

The entire notion of ultraposets, or more generally ultracategories, is an axiomatisation of ultraproduct structures. Such duality results again shows the significance of \L os Ultraproduct Theorem in logic: From the ultrastructure on the space of points, we can reconstruct the syntax of the theory upto equivalence.

It is also clear from this paper that the notion of \emph{left} and \emph{right} ultrafunctors, not just ultrafunctors, are absolutely essential to our development, and they indeed incorporates many different notions in different contexts, as we have partially seen here. This is the reason why we think Lurie's approach~\cite{lurie2018ultracategories} to ultracategories is very promising.

As we have mentioned in the introduction, we consider this paper as a first step into a more detailed study of ultraposets and ultracategries in general, with potential applications to logic. We end the paper by discussing several future directions naturally arise in this paper.

\subsection{Ultraposets and Ordered Topological Structures}

We hope that this paper has made the connection between ultraposets and ordered topological structures clear. We have seen from Lemma~\Ref{lem:compord} that the category $\mb{CompOrd}$ of compact ordered spaces fully faithfully embeds into $\UltPos$, and this embedding brings a lot of examples we encounter in practice to the world of ultraposets.

In particular, Section~\Ref{sec:sep} has mentioned that all Priestley spaces are of this form. However, in this case, Theorem~\Ref{thm:cloccluteq} makes it clear that considering \emph{left} ultrafunctors between these spaces is still interesting. This aspect cannot be seen from the perspective of ordered topological structure alone, because as we have seen, left ultrafunctors might even not be continuous for these topological spaces. This makes it interesting to see how the theory of other ordered topological structures can be studied by viewing them as ultraposets, and whether there's further connections with other fields, like locale theory.

\subsection{Ultraposets and Domain Theory}

As we have seen, Lemma~\Ref{lem:leftscott} also suggests there is a fully faithful embedding of $\CPos$, the (2-)category of complete lattices with Scott-continuous maps between them, into $\UltLPos$.

Of particular importance to Domain theory is the full subcategory of \emph{continuous lattices} of $\CPos$. In particular, every Scott-continuous function between continuous lattices will be a \emph{left ultrafunctor}, and if the function furthermore preserves arbitrary meets, it will in fact be an \emph{ultrafunctor}. These notions of morphisms between continuous lattices are widely used in domain theory (cf.~\cite{gierz2003continuous}) to construct models of various type theories, and can indeed be subsumed in the theory of ultraposets. It will be interesting to see whether further domain theoretic properties can be expressed in the language of ultraposets, and how that could be useful for domain theory.

\subsection{Ultraposets and Locales}

In previous sections, we have constructed two functors
\[ \Omega_u : \UltPos \to \DL\op, \quad \Omega : \UltLPos \to \Loc. \]
We have shown that the former in fact has a right adjoint $\Mod$, and this identifies $\DL\op$ as a reflexive subcategory of $\UltPos$. However, the latter $\Omega$ does not have a right adjoint immediately, since the poset of points $\pt(L)$ of an arbitrary locale may not be closed in $[L,\mb 2]$ equipped with the canonical ultrastructure.

However, we have mentioned in Remark~\Ref{rem:betaloc} that \cite{di2022geometry} has identified a class of locales whose poset of points indeed inherits a canonical ultrastructure. These are locales that are right Kan injective w.r.t. embeddings of spaces $X \hook \beta X$ for any set $X$, viewed as a localic map from the discrete space $X$ to the compact Hausdorff space $\beta X$. Let $\beta\Loc$ be the category of these locales. Is $\Omega$ part of an adjoint between $\UltLPos$ and $\beta\Loc$? We leave this for future work.

\subsection{Theoretic Development of Ultraposets}

The most important reason we have adopted Lurie's approach to ultracategories, and ultraposets in particular, is that it naturally leads to the notion of \emph{left} and \emph{right} ultrafunctors, which are more fundamental than the notion of ultracategories.

Although we have studied ultraposets and (left/right) ultrafunctors in general in this paper, to a large extent we have mainly focused on zero-dimensional ones. However, these are quite restrictive in the class of all ultraposets, because each zero-dimensional ultraposet comes from a compact ordered space. 

However, quite a lot of examples coming from different fields naturally embeds in one of the categories of ultraposets $\UltPos$, $\UltLPos$ or $\UltRPos$. In general they will be quite different in nature. A more systematic study of ultraposets themselves has the potential of providing a unifying treatment of these different topics, and we leave this for future development.

\bibliography{mybib}{}
\bibliographystyle{apalike}

\appendix

\section{Properties on Ultrafilters and Lattices}\label{app}

In this appendix, we record some basic properties about ultrafilters and lattices that are used in this paper. The following lemma is crucial:

\begin{lemma}\label{lem:dismax}
  Given any distributive lattice $D$, let $I$ be an ideal and let $F$ be a filter in $D$. If $I$ and $F$ are disjoint, then $I$ can be extended to a maximal ideal that is disjoint from $F$. Similarly, $F$ can be extended to a maximal filter that is disjoint from $I$. These maximal ideals and filters will be prime.
\end{lemma}
\begin{proof}
  We only prove the case for extensions of ideals; the case for filters is completely similar. The existence of such maximal ideals disjoint from $F$ is implied by Zorn's Lemma. For any such maximal ideal $I$ disjoint from a filter $F$, let $p,q\in D$ that $p \wedge q$ in $I$. Suppose $p\not\in I$, then the ideal $I\pair p$ obtained by joining $p$ to $I$ would property extend $I$, thus by maximality of $I$ it would intersect $F$. Concretely, this implies there exists $i\in I$ that $i \vee p \in F$. Similarly, if $q \not\in I$, then there exists $j\in I$ that $j \vee q \in F$. Now suppose both $p,q\not\in I$, then we have
  \[ (i\vee p) \wedge (j \vee q) \in F. \]
  On the other hand, we have
  \[ (i \vee p) \wedge (j \vee q) = ((i\vee p) \wedge j) \vee ((i\vee p) \wedge q). \]
  Evidently, $(i\vee p) \wedge j \in I$. We also know that
  \[ (i\vee p) \wedge q = (i \wedge q) \vee (p \wedge q) \in I, \]
  because by assumption $p \wedge q \in I$. This actually implies that
  \[ (i\vee p) \wedge (j \vee q) \in I, \]
  contradicting the assumption that $I$ and $F$ are disjoint. Thus, at least one of $p,q$ is in $I$, thus $I$ is a prime ideal.
\end{proof}

\subsection{Ultrafilters}

For any set $X$, an ultrafilter $\mu$ on $X$ is a family of subsets of $X$, such that it is upward closed, closed under finite intersection, and for any subset $A$ of $X$, either $A\in\mu$ or $\ov A \in \mu$. Equivalently, $\mu$ is a \emph{maximal filter} on $\wp(X)$.

Let $X$ be a set, we say a family $\set{K_i}$ of subsets of $X$ has the \emph{finite intersection property} if any finite subfamily of $\set{K_i}$ has non-empty intersection.

\begin{lemma} \label{lem:maximalideal}
  For any set $X$ and any family $\set{K_i}$ of $X$ that has the finite intersection property, there exists an ultrafilter $\mu$ on $X$ that extends this family, in the sese 
\end{lemma}
\begin{proof}
  Any such family can be extended to a filter, by closing them under finite intersection. Then by Lemma~\Ref{lem:dismax}, this filter can be extended to a maximal filter on $X$, which is exactly an ultrafilter. 
\end{proof}

\subsection{Distributive Lattices and their Models}

Models of a distributive lattice corresponds to its prime ideals:

\begin{lemma}\label{lem:modelprim}
  A monotone map $x : D \to \mb 2$ is a morphism between distributive lattices iff $x\inv(0)$ is a prime ideal in $D$.
\end{lemma}
\begin{proof}
  Suppose $x : D \to \mb 2$ is a morphism of distributive lattices. For any $p,q\in D$, if $p \wedge q \in x\inv(0)$, then $x(p) \wedge x(q) = x(p \wedge q) = 0$. It follows that at least one of $x(p),x(q)$ equals to 0, hence $p \in x\inv(0)$ or $q\in x\inv(0)$.

  On the other hand, if $x\inv(0)$ is a prime ideal, then for any $p,q\in D$, if $p,q\in x\inv(0)$, then both $p \wedge q \in x\inv(0)$ and $p \vee q \in x\inv(0)$. If $p\in x\inv(0)$ and $q \not\in x\inv(0)$, then $p \wedge q \in x\inv(0)$, while $p \vee q\not\in x\inv (0)$, because $x\inv(0)$ is downward closed. If $p,q\not\in x\inv(0)$, then both $p\wedge q\not\in x\inv(0)$ because $x\inv(0)$ is prime, and $p\vee q\not\in x\inv(0)$. In each case, $x$ preserves the lattice structure. 
\end{proof}

The following theorem is equivalent to the completeness theorem of propositional logic:

\begin{lemma} \label{lem:idealintersecprime}
  Let $I$ be an ideal of a distributive lattice $D$, then $I$ is the intersection of all prime ideals extends $I$.
\end{lemma}
\begin{proof}
  For any $p\not\in I$, consider the principle filter $\dv p$ generated by $p$, which by assumption is disjoint from $I$. By Lemma~\Ref{lem:dismax}, we can extend $I$ to a maximal ideal $x$ disjoint from $\dv p$, and this ideal $x$ would be prime. Thus, for any $p\not\in I$, there exists some prime ideal $x$ containing $I$ that does not contain $p$, hence $I = \bigcap_{I\subseteq x} x$. 
\end{proof}

\begin{corollary}[Completeness of Propositional Logic]\label{cor:completeness}
  For any distributive lattice $D$, if $p \not\le q$, then there exists some prime ideal $x$ of $D$, such that $q\in x$ and $p \not\in x$.
\end{corollary}
\begin{proof}
  Consider the ideal $\cv q$ in $D$. Since $p \not\le q$, then $p \not\in \cv q$. By Lemma~\Ref{lem:modelprim}, there must exist some prime ideal $x$ extending $\cv q$, such that $p \not\in x$.
\end{proof}

\end{document}